\documentclass[12pt]{article}

\usepackage{amssymb,amsmath,amsthm, amstext}
\usepackage{latexsym}
\usepackage{natbib}
\usepackage{amsbsy}
\usepackage{graphicx}
\usepackage{float}
\usepackage{subfloat}
\usepackage{caption}
\usepackage{subcaption}
\usepackage{epsfig}
\usepackage{epstopdf}  
\usepackage{color}
\usepackage[font= {small,it} ]{caption}
\usepackage{url}
\usepackage[titletoc]{appendix}
\usepackage{xcolor}  
\usepackage{hyperref}
\usepackage[margin= 1in]{geometry}
\usepackage{appendix}
\usepackage{bm}
\usepackage{multirow}

\usepackage{algorithmic}

\usepackage[ruled,vlined]{algorithm2e}

\usepackage{tikz-cd}

\newcommand*\samethanks[1][\value{footnote}]{\footnotemark[#1]}

\newtheorem{theorem}{Theorem}[section]
\newtheorem{lemma}[theorem]{Lemma} 
\newtheorem{proposition}[theorem]{Proposition}

\newtheorem{definition}[theorem]{Definition}

\begin{document}

\title{A Kernel-Based Approach for Gaussian Process Modeling with Functional Information\thanks{{This work was funded by the National Science Foundation under grant DMS 2210686}}}
\author{D. Andrew Brown\thanks{Corresponding author, School of Mathematical and Statistical Sciences, Clemson University, Clemson, SC 29634, USA, ab7@clemson.edu; Order of authorship is alphabetical} \and Peter Kiessler\thanks{School of Mathematical and Statistical Sciences, Clemson University, Clemson, SC 29634, USA} \and John Nicholson\samethanks}

\maketitle

\begin{abstract}
Gaussian processes (GPs) are ubiquitous tools for modeling and predicting continuous processes in physical and engineering sciences. This is partly due to the fact that one may employ a Gaussian process as an interpolator while facilitating straightforward uncertainty quantification at other locations. In addition to training data, it is sometimes the case that available information is not in the form of a finite collection of points. For example, boundary value problems contain information on the boundary of a domain, or underlying physics lead to known behavior on an entire uncountable subset of the domain of interest. While an approximation to such known information may be obtained via pseudo-training points in the known subset, such a procedure is {\em ad hoc} with little guidance on the number of points to use, nor the behavior as the number of pseudo-observations grows large. We propose and construct Gaussian processes that unify, via reproducing kernel Hilbert space, the typical finite training data case with the case of having uncountable information by exploiting the equivalence of conditional expectation and orthogonal projections in Hilbert space. We show existence of the proposed process and establish that it is the limit of a conventional GP conditioned on an increasing number of training points. We illustrate the flexibility and advantages of our proposed approach via numerical experiments.

\end{abstract}

\section{Introduction}\label{sec:intro}
Gaussian processes (GPs) \cite{rasmussen06} are popular tools among scientists and engineers for modeling complex physical processes because of their flexibility, simplicity, and closed-form quantification of uncertainty. They are commonly employed as surrogate models to be run in place of computationally expensive computer models \cite{gramacy2020surrogates}. (Polynomial chaos expansions and neural networks are also often used, but these have been shown to perform similar to or even worse than GPs \cite{OwenEtAl17, myrenGPvsNN}.) As Gaussian processes have become more popular in practice, there has arisen a demand to modify GPs to possess certain characteristics; e.g., to honor known physics \cite{RaissiEtAl17}, shape constraints \cite{WangBerger16}, orthogonality \cite{Plum18}, or boundary constraints \cite{ding2019bdrygp}. Such modifications are useful for improving the interpolation performance of the GPs as well as mitigating identifiability issues that arise in, e.g., calibration of computer models \cite{koh01, BryOH14}. A review of constrained GPs may be found in \cite{swiler}.

In differential equations, boundary constraints on the values of the solution are called Dirchlet boundary conditions (as opposed to, e.g., Neumann boundary conditions which specify values of the derivatives). This is a common setting for modeling GPs. In a more general scenario, however, one may simply have knowledge of a process on a subset of the domain. This does not always fit under the umbrella of ``boundary conditions," as the knowledge of the process may not be on the boundary. In this paper, we propose and explore theoretically and numerically a novel class of Gaussian processes which have known, fixed values on an arbitrary compact subset of the domain.

As motivation, consider the following scenario that arises in materials science. Finite element models can be used to predict the strength of composite materials consisting of a polymer matrix and a filler material made of embedded spherical particles \cite{arp21}. There are seven parameters contributing to variations in strength, six of which determine properties of the filler and interactions between the filler and the matrix. The finite element model code is too computationally expensive to run directly, so a GP can serve as a surrogate for the model output. When there is no filler in the material, though, the strength of the composite is simply the strength of the polymer, which is entirely known and controllable. In other words, the strength of the composite is known on an uncountable, six-dimensional subset of the seven-dimensional domain. A temptation in this situation is to include a finite number of ``pseudo observations" along the known subset as additional training data for the GP, since these training points are available at no additional cost. However, this still raises the question of how many of these points to use and how they should be distributed throughout this subset. Indeed, any finite number of points does not completely capture all the available information. Our goal in the present work is to construct and study a method for more fully incorporating such {\em a priori} knowledge into GPs to capture such information in a more principled way. 

The idea behind our proposed construction is that fixing the value of a Gaussian process at certain points can be framed as finding the conditional distribution. For Gaussian distributions, conditioning on a finite number of points is well-known and follows from standard multivariate normal theory. Conditioning on uncountable subsets, however, is not as straightforward. Our approach is to view conditional expectation as an orthogonal projection so that determining the conditional distribution reduces to explicitly identifying the form of the projection, which we are able to do. 

\paragraph{Related Work} There exist in the literature several proposed approaches for incorporating boundary constraints into GPs, which is a special case of the problem considered in this work. \cite{Solin_2019} suggest modifying an analytic stationary covariance function by approximation with a collection of basis functions obtained via spectral decomposition of the homogenous Laplace equation, basis functions that vanish on the boundary of the domain. \cite{hegermann} use pushforward GP mappings of the form $\rho \mathbb{X}$, where $\rho:\mathbb{R}^d \to [0,1]$. The author suggests choosing $\rho$ so that $\rho \equiv 0$ on the boundary as a means of satisfying the constraint. In preceding work, \cite{tan} developed an explicit construction following the same reasoning as that of \cite{hegermann}, and developed a mean function which permits nonzero constant boundary conditions. \cite{ding2019bdrygp} defined a boundary-valued GP with a covariance function that vanishes on all or part of the boundary, yielding the known values contained in the mean function. Though these methods have proven reasonable and effective under certain circumstances, none are able to handle more general domain constraints. 

Similar to the work in this paper, \cite{Owhadi2018} also exploit the equivalence between conditional expectation and orthogonal projection to establish the existence of Gaussian measures on orthogonally-decomposed Hilbert spaces. Specifically, \cite{Owhadi2018} define a sequence of projection operators associated with conditional mean and covariance functions, and show weak operator convergence to the mean and covariance operators conditioned on infinite-dimensional space. While existence of conditional Gaussian measures with these mean and covariance functions is established, the authors do not consider the sequence of stochastic (Gaussian) processes associated with the sequence of finite-dimensional operators and the process associated with the limiting operator. Recently, \cite{steinwart2024} extended the ideas of \cite{Owhadi2018} to more general Banach space by using filtering sequences in place of projections to arrive at the same convergence results. Further, \cite{steinwart2024} showed weak convergence of the associated stochastic processes to the process associated with the limiting operators. To the best of our knowledge, \cite{steinwart2024} is the only work other than ours presented here to show weak convergence of the conditional Gaussian processes. While both \cite{Owhadi2018} and \cite{steinwart2024} provide important theoretical results, they are quite limited in practical guidance on the application of such results to problems often encountered in uncertainty quantification (UQ). 

\paragraph{Contributions} In this work we likewise consider and establish Gaussian processes conditioned on an increasingly dense subset of an arbitrary compact set contained in the domain of interest. This facilitates convergence of the finite-dimensional distributions to that of a GP conditioned on values on the entire compact set. Similar to \cite{steinwart2024}, we go beyond finite-dimensional convergence and show weak convergence of the sequence of stochastic processes themselves to the limiting GP. An important distinction between the aforementioned earlier works and our own, though, is that we construct the reproducing kernel Hilbert space (RKHS), $\mathcal{H}(T)$, in terms of the covariance function.  In our setting, as well as those of \cite{Owhadi2018, dacosta2024, steinwart2024}, the RKHS exists and is exactly the Cameron-Martin space as defined by \cite[p. 44]{Bog1998}. The RKHS setting deserves special attention for several reasons.
First, it allows us to prove very similar results to earlier work \cite{Owhadi2018, dacosta2024, steinwart2024} but in a more straightforward (and faster) manner. For instance, we are able to almost immediately prove weak convergence of the sequence of GPs as a direct consequence of the Kolmogorov-Chentsov conditions after establishing the existence of continuous versions and tightness. Second, compared to \cite{Owhadi2018}, we can derive stronger results; e.g., we use the RKHS properties to show strong operator convergence of the projections. Third, and perhaps most importantly, the RKHS assumption allows us to propose and illustrate practical application of these results. In our case, we leverage the algorithm of \cite{oya2009} for approximating RKHS inner products to propose and illustrate guidance on numerical implementation, and how our approach can improve on alternatives.

{The reason the calculations simplify in our setting is as follows. Consider a Gaussian process on a compact set $T \subset \mathbb{R}^d$, $\{X_t : t \in T\},$ with covariance kernel $k$. If $k$ is smooth enough (in a sense to be made more precise in the sequel), there is a version of $\{X_t : t \in T\}$ that is continuous (See Lemma \ref{lem:continuous} in Subsection \ref{sec:prob} below.) In this case, the distribution of the GP, $\Gamma$, can be considered as a probability measure on $C(T)$. But $C(T)$ can be continuously embedded as a dense subset of $L^2(T)$. Consequently, $\Gamma$ can also be considered as a measure on $L^2(T)$. Since $L^2(T)$ is a Hilbert space, we can identify its dual space with $L^2(T)$. Thus, each $f \in L^2(T)$ can be viewed as a Gaussian random variable, and for $f, g \in L^2(T)$,
\begin{equation}\label{eqn:l2cov}
    \mbox{Cov}(f, g) = \iint_{T \times T} f(s)k(s,t)g(t)dsdt.
\end{equation}
This means that, for each $f \in L^2(T)$, $\tilde{f} := \int_T k(\cdot,t)f(t)dt$ belongs to $\mathcal{H}(T)$ and the collection of such $\tilde{f}$ is dense in $\mathcal{H}(T)$. Moreover, the inner product of $\tilde{f}$ and $\tilde{g}$ in $\mathcal{H}(T)$ is given by \eqref{eqn:l2cov}. Thus, calculations in the Hilbert space $\mathcal{H}(T)$ can be facilitated by the integral operator $K$ of $L^2(T)$ defined by $Kf(s) = \int_T k(s,t)f(t)dt = \tilde{f}(s)$. The operator $K$ has many useful properties that we can exploit for our analysis in this work.}

To fix ideas, consider a Gaussian field $\mathbb{X}^0=\{X^0_s: \, s \in T\}, ~T \subset \mathbb{R}^d$, with mean function $\mu$ and covariance kernel $k$. For $n$ discrete points $t_1, \ldots, t_n \in T$, it is well-known that the process $\mathbb{X}^n=\{X^n_s: \, s \in T\}$ , where $X^n_s:=X^0_s|(X_{t_1}=x_{t_1},...,X_{t_n}=x_{t_n})$,
is also a Gaussian process with mean function $\mu$
\begin{equation}\label{eq:finiteMean}
    \mu_0(\cdot) =\mu(\cdot)+k(\cdot,\mathbf{t}) k(\mathbf{t},\mathbf{t})^{-1}(\mathbf{x}-\mu(\mathbf{t})),
\end{equation}
and covariance kernel
\begin{equation}\label{eq:finiteKernel}
    k_0(\cdot,\cdot) =k(\cdot,\cdot)-k(\cdot,\mathbf{t}) k(\mathbf{t},\mathbf{t})^{-1} k(\mathbf{t},\cdot),
\end{equation}
where $\mathbf{t} = (t_1, \ldots, t_n)^\top$ and $\mathbf{x} = (x_{t_1}, \ldots, x_{t_n})^\top$. This can be derived by projecting the (unconditional) mean function $\mu(\cdot)$ and covariance kernel $k(\cdot, \cdot)$ onto the function space associated with $T_0 = \{t_1, \ldots, t_n\} \subset T$. 

In the finite dimensional case, projections typically can be computed explicitly using elementary linear algebra \cite{stapleton}. For infinite dimensional function spaces, our approach in this work is similarly to associate to the distribution of a Gaussian process $\mathbb{X}^0$ conditional on $\mathbb{X}^0|_{T_0}=g_0$ an orthogonal projection from one function space to another, where $T_0 \subset T$ is the set on the which the values of the GP are known exactly. We rigorously describe the projection operator and use it to find the conditional distribution. In the process, we show that our approach unifies conditioning on a finite set of points with that on an uncountable, compact subset of the input space. We find the conditional mean and covariance functions and show that the associated GP does, in fact, exist. Further, we formally establish that the resulting GP, which we term projected kernel Gaussian process (pkGP), is the limit of GPs conditioned on an increasing number of finite points in the known subset, following our intuition.

This paper is organized as follows: Section \ref{sec:prelim} reviews the pertinent ideas from the theory of reproducing kernel Hilbert spaces (RKHSs) \cite{paulsen} and derives the well-known finite dimensional conditional distribution from the RKHS perspective. Section \ref{sec:general} presents our results in the general setting, including existence and weak convergence of the associated GPs. This section also briefly discusses considerations associated with adding a nugget to the covariance function, as commonly done to improve the condition numbers of matrices associated with certain kernels. Section \ref{sec:rrAndKrig} discusses computational implementation of the RKHS inner products, including an illustration. Section \ref{sec:simulations} contains the results of numerical experiments in which we interpolate several different test functions with different types of domain constraints. The paper concludes with some final remarks in Section \ref{sec:final}. Throughout this work, we draw on several fundamental results from probability, functional analysis, and RKHS theory that can be found in, e.g., \cite{kall}, \cite{lax}, and \cite{paulsen}, respectively.

\section{Preliminaries}
\label{sec:prelim}

{We begin this section with a formal definition of a Gaussian process.}
\begin{definition}
    {For a measurable input space $(T, \mathcal{G})$ and a probability space $(\Omega, \mathcal{F}, P)$, a Gaussian process is a stochastic process $\mathbb{X} = \{X(t, \cdot) : t \in T\}$ such that for a fixed $\omega \in \Omega$, $X(\cdot, \omega): T \rightarrow \mathbb{R}$ is a $\mathcal{G}-$measurable function, and for any $\boldsymbol{t} = (t_1, \ldots, t_N)^\top \in T^N$, $X(\boldsymbol{t}, \cdot): \Omega \rightarrow \mathbb{R}^N$ follows a multivariate Gaussian distribution. For simplicity, we denote this as $\mathbb{X} = \{X_t : t \in T\},$ where $X_t := X(t, \cdot)$.}
\end{definition}
In our setting, Gaussian processes are typically used to learn continuous, often differentiable functions via conditioning their distributions on known locations and function values, determined by \eqref{eq:finiteMean} and \eqref{eq:finiteKernel}. Likewise, in our work we make the often reasonable assumption that the target function is continuous. However, orthogonal projections as mentioned in Section \ref{sec:intro} are not permissible within the space of continuous functions, $C(T)$, since $C(T)$ {is not a Hilbert space}. On the other hand, {by assuming as we do here that the input domain is compact and the kernel is continuous}, reproducing kernel Hilbert spaces (RKHS) \cite{paulsen} are subsets of $C(T)$ containing functions that, under modest conditions, can serve as approximations to other functions to an arbitrary degree of precision. As Hilbert spaces with associated inner products, orthogonal projections can be defined on them. Therefore, the theory developed in this paper will use analytical and probabilistic properties of RKHSs.

\subsection{Definition and Overview} \label{sec:rkhs_intro}

Construction of a {Gaussian} conditional distribution revolves around an appropriate covariance function, which for the case of Gaussian processes (GPs) will be studied as an element of a function space. In this section we briefly review RKHSs, integral operators, and how orthogonal projection in Hilbert space leads to the well-known GPs conditional on a finite number of observed values. 

{Take $T$ to be a compact subset of $\mathbb{R}^d$.} Let $k: T\times T \rightarrow \mathbb{R}^+$ denote the covariance function of a Gaussian process.  As such, it is symmetric in its arguments and positive definite. We assume further that it is continuous.  Let $K$ denote the integral operator in $L^2(T)$ associated with the kernel $k$, defined by 
\begin{equation}
Kx(t)= \int_T k(s,t) x(s) ds. \label{eq:integral_operator}
\end{equation}
We denote  the range of $K$ as $R(K)$ and define {$\langle  \cdot \, , \cdot \rangle_T$} to be the standard inner product on $L^2${; i.e., {$\langle f, g \rangle_T = \int_T f(s)g(s) ds$}}.

For $t \in T$, define $\varrho_t: f \mapsto f(t)$ to be the {evaluation} functional. These are commonly seen defined on $(C(T),||\cdot||_\infty)$ where $||\cdot||_\infty$ denotes the supremum norm. As elements of the dual space, the evaluation functionals correspond to Dirac measures. The motivation behind RKHS is to construct a Hilbert space so that each evaluation functional is bounded and thus identifies uniquely with an element of the space itself. Thus, a RKHS on $T$, $(\mathcal{H}(T),\langle \cdot,\cdot \rangle_{\mathcal{H}(T)}),$ is defined to be the collection of functions such that the evaluation functionals are bounded. 

A kernel $k$ defined on $T \times T$ has the reproducing property on $\mathcal{H}(T)$ if the representation of $\varrho_t$ in $\mathcal{H}(T)$ is {the kernel section} $k_t:=k(\cdot,t)$ for each $t \in T$. It follows that the inner product $\langle \cdot, \cdot \rangle_{\mathcal{H}(T)}$ satisfies $f(t)=\langle f, k_t\rangle_{\mathcal{H}(T)}$, for any $f \in \mathcal{H}(T)$ and $t \in T$. By the Moore-Aronszajn  Theorem, each RKHS is identified uniquely with a kernel \cite[Theorem 2.14]{paulsen}. The {RKHS associated with $k$} is constructed by closing the span of the functionals $\{k_t\}_{t \in T}$ under $||\cdot||_{\mathcal{H}(T)}$, implying of course that $\{k_t\}_{t \in T} \subset \mathcal{H}(T)$. In addition, the norm of $k_t$ can be calculated explicitly by $||k_t||_{\mathcal{H}(T)}=\langle k_t,k_t\rangle^{1/2}_{\mathcal{H}(T)}=k(t,t)^{1/2}$.
Furthermore, for $s,t \in T$,
\begin{align*}
    ||k_s-k_t||^2_{\mathcal{H}(T)} = \langle k_s - k_t,k_s-k_t \rangle_{\mathcal{H}(T)} = k(s,s)-k(s,t)-k(t,s)+k(t,t).
\end{align*}
{This establishes the following Lemma}:
\begin{lemma}\label{lemma:kHolder}
If $k$ is $\gamma$-H\"older continuous, then $||k_s-k_t||^2_{\mathcal{H}(T)} \le B |s-t|^{\gamma}$, for some constant $B>0$.
\end{lemma} 
This fact plays an important role in Section \ref{sec:prob}, where we show weak convergence of Gaussian processes to a limit.

Mercer's theorem \cite[p. 343]{lax} plays a fundamental role in the theory of RKHS. It states that if $k$ is a continuous kernel, then for any $s,t \in T$, 
$$
k(s,t)=\sum_{n=1}^\infty \lambda_n e_n(s)e_n(t),
$$
{where $\{(\lambda_n, e_n)\}_{n=1}^{\infty}$ are the pairs of eigenvalues and orthonormal eigenfunctions associated with $K$, and} the series converges absolutely and uniformly. In addition, it can be shown that for $f,g \in \mathcal{H}(T)$,
\begin{equation}\label{eq:rkhs_ip}
\langle f,g \rangle_{\mathcal{H}(T)}=\sum_{n=1}^\infty \frac{{\langle f,e_n \rangle_T\langle g,e_n \rangle_T}}{\lambda_n}, 
\end{equation}
and thus any $f \in \mathcal{H}(T)$ must satisfy $\sum_{n=1}^\infty \frac{{\langle f,e_n \rangle_T^2}}{\lambda_n} < \infty$. We can generalize this to say that $\mathcal{H}(T)=\{\sum_{n=1}^\infty a_n e_n: \sum_{n=1}^\infty \frac{|a_n|^2}{\lambda_n} < \infty \}$.

Consider the square root operator $K^{1/2}$ of the integral operator $K$. {Since $k(\cdot, \cdot)$ is continuous, $K : L^2(T) \rightarrow L^2(T)$ and thus $K^{1/2}:L^2(T) \rightarrow L^2(T)$ are compact \cite[Theorem 4.6.2]{hsingEubank15}. Further, we assume $k(\cdot, \cdot)$ is symmetric in its arguments, whence $K^{1/2}$ is also self-adjoint. The square root operator can be expressed as \cite[p. 100]{hsingEubank15}}
$$
K^{1/2} x = \sum_{n=1}^\infty \lambda_n^{1/2} {\langle x, e_n \rangle_T} e_n, {~~\forall x \in L^2(T)}.
$$
For $x \in L^2(T),~||K^{1/2}x||_{\mathcal{H}(T)}^2  = \langle K^{1/2}x,K^{1/2}x \rangle_{\mathcal{H}(T)}
 = \sum_{n=1}^\infty {\langle x,e_n \rangle_T^2} \le ||x||_{L^2}^2$, by Bessel's inequality. In particular, if $K$ has a trivial nullspace, the eigenvectors $\{e_n\}$ form {an orthonormal} basis of $L^2(T)$, which allows us to substitute the inequality with an equality. If this is the case, $K^{1/2}$ is an isometric isomorphism between $L^2(T)$ and $\mathcal{H}(T)$. Hence, $K^{-1/2}$ exists and is bounded, and for $f,g \in \mathcal{H}(T),$
\begin{equation}\label{eq:rootInnerProduct}
\langle f,g \rangle_{\mathcal{H}(T)}= {\langle K^{-1/2}f,K^{-1/2}g \rangle_T}.
\end{equation}
Note that $K^{1/2}: L^2(T) \rightarrow \mathcal{H}(T)$ is bijective due to the restriction from $L^2(T)$ to $\mathcal{H}(T)$, which loses compactness of the operator but allows for the existence of the inverse $K^{-1/2}$.

The projection occurs in both the mean and the covariance, meaning that the mean function should be an element of the RKHS. If the mean function is zero, this is trivially the case. Otherwise, it is difficult to check if a function is an element of $\mathcal{H}(T)$. For example, it has been shown that the RKHS associated with the square exponential kernel, $k(s,t)=\exp \{-|s-t|^2\}$, does not contain any {non-zero} constant functions or polynomials \cite{ha}. When the mean function is not an element of the RKHS, it is important that it can be well approximated by an element of the RKHS. 

A key but not restrictive assumption that we make throughout this paper is that the kernel used is universal \cite{micchelli}. A kernel is said to be universal if for any compact subset $Z$ of the input space, the RKHS it generates is dense in the continuous functions on $Z$ under the supremum norm. This class of kernels includes all power exponential kernels of the form $k(s,t) = \exp\{\ell |s - t| ^ p\}, \, \ell, p >0$, as well as the M\'atern and rational quadratic kernels.

\subsection{Finite Case}\label{sec:ortho}

Orthogonal projections in an RKHS are not as simple to visualize as they are in finite dimensional or $L^2$ spaces. However, they have important properties for our purposes regarding the connection to their respective generating kernels.

Suppose that $P:\mathcal{H}(T) \to \mathcal{H}_0$ is {the orthogonal projection into a {closed} subspace $\mathcal{H}_0 \subset \mathcal{H}(T)$}, keeping in mind that we are interested in subspaces of the form $\mathcal{H}_0 = \{f \in \mathcal{H}(T) : f|_{T_0} \equiv 0\}$, where $T_0 \subset T$. {For $s, t \in T$, we can define $Pk$ by its action on the kernel sections $k_s$ and $k_t$.} By properties of orthogonal projections, we have that $Pk(s, t) = \langle Pk_s, k_t \rangle_{\mathcal{H}(T)} = \langle Pk_s, Pk_t \rangle_{\mathcal{H}(T)} = \langle Pk_s, Pk_t \rangle_{\mathcal{H}_0}$. This leads us to an important result regarding how $\mathcal{H}_0$ is generated, the proof of which can be found in \cite[Theorem 2.5]{paulsen}.

\begin{proposition}
    $\mathcal{H}_0$ is a RKHS with reproducing kernel $k_0(s, t) = Pk(s, t)$.
\end{proposition}

Now take {$\{t_1, \ldots, t_n\} = T_0$ and} $\mathcal{H}_0 = \{ f\in \mathcal{H}(T): f(t_i)=0, ~i=1,\ldots,n\}$.
As we will show in Proposition \ref{prop: perp}, $\mathcal{H}_0^\perp = {\overline{\mbox{Span}}}(\{k_{t_1},\ldots, k_{t_n}\}) ~{= {\mbox{Span}}(\{k_{t_1},\ldots, k_{t_n}\})}$. 
As closed subspaces of $\mathcal{H}(T)$, both $\mathcal{H}_0$ and $\mathcal{H}_0^\perp$ are RKHSs.
It is easier to find the kernel that generates $\mathcal{H}_0^\perp$, so we do that en route to finding the kernel that generates $\mathcal{H}_0$.  Toward this end, let $Q$ be the orthogonal projection onto $\mathcal{H}_0^\perp$. Then, for $f \in {\mathcal{H}(T)}$, 
\begin{equation} \label{eq:finite_proj}
    Qf(\cdot) = \sum_{i = 1} ^ n a_i k_{t_i}(\cdot), \, a_i \in \mathbb{R}, \, i = 1,...,n. 
\end{equation}
Observing that $Qf(t_i) = \langle Qf, k_{t_i} \rangle_{\mathcal{H}(T)} = \langle f, Qk_{t_i} \rangle_{\mathcal{H}(T)} = \langle f, k_{t_i} \rangle_{\mathcal{H}(T)} = f(t_i),$ it follows that $Qf$ is an interpolation of $f$ at the points $\{t_i\}_{i=1}^n$. Defining $k(\mathbf{t}, \mathbf{t}) = (k(t_i, t_j))_{i,j=1}^n$, {$\mathbf{a} = (a_1, \ldots, a_n)^\top$}, and $f(\mathbf{t}) = (f(t_1), ..., f(t_n))^\top$,
it follows that 
\begin{equation}\label{eq:finite_proj_coef}
    \mathbf{a} = k(\mathbf{t}, \mathbf{t})^{-1} f(\mathbf{t}).
\end{equation}

\noindent Choosing $f = k_{s_1}$ for $s_1 \in T$, and using \eqref{eq:finite_proj} and \eqref{eq:finite_proj_coef}, we have that 
\begin{equation} \label{eq:finite_proj_kernel}
    Qk(s_1, s_2) = Qk_{s_1}(s_2) = k(s_1, \mathbf{t}) k(\mathbf{t}, \mathbf{t})^{-1} k(\mathbf{t}, s_2).
\end{equation} 
One may recognize that the righthand side of this equation appears in \eqref{eq:finiteKernel}. 

Using the decomposition $\mathcal{H}(T) = \mathcal{H}_0 \oplus \mathcal{H}_0^\perp$, we can say the following about the kernels of $\mathcal{H}(T), \, \mathcal{H}_0, \, \mathcal{H}_0^\perp$. The proof can be found in \cite[Corollary 5.5]{paulsen}:
\begin{proposition}
    {Let $k_0$ be the kernel which generates $\mathcal{H}_0$ and $k_\perp$ the kernel which generates $\mathcal{H}_0^\perp$. Then, $k = k_0 + k_\perp$ and therefore $k_0 = k - k_\perp$}.
\end{proposition}
{Hence, by definition of $Q$}, we have
\begin{equation} \label{eq: cov_projection}
    k_0 = k - Qk.
\end{equation}


Lastly, let $g ~{\in \mathcal{H}(T)}$ represent the function upon which we want the GP to be fixed at $\{t_1, ..., t_n\}$. Then it is necessary for the GP mean function $\mu$ to be restricted to a conditional mean $\mu_0$ satisfying $\mu_0(t_i) = g(t_i),\, i = 1, \ldots, n$; i.e., $\mu_0 - g \in \mathcal{H}_0$. {Likewise, $\mu - Q\mu \in \mathcal{H}_0$ and $\mu(\mathbf{t}) - Q\mu(\mathbf{t}) = \mu_0(\mathbf{t}) - g(\mathbf{t})$. Assuming $g \in \mathcal{H}_0^\perp$ so that $g = Qg$ (a reasonable assumption given the definition of {$\mathcal{H}_0^\perp$}), we have $\mu_0(\mathbf{t}) = \mu(\mathbf{t}) + Q(g(\mathbf{t}) - \mu(\mathbf{t}))$. Hence,}
\begin{equation} \label{eq: mean_projection}
    \mu_0 = \mu + Q(g - \mu),
\end{equation}
which is analogous to \eqref{eq:finiteMean}. Thus, given a GP $\mathbb{X}$ on $T$ with mean $\mu$ and covariance $k$, using orthogonal projections on RKHS we are able to modify $\mathbb{X}$ so that $X(t_i) = g(t_i), \, i = 1,...,n$. The resulting process is determined by the mean function $\mu + Q(g - \mu)$ and covariance $k - Qk$. This formulation will be shown in the sequel to remain true when considering more general subsets of $T$. In particular, Section \ref{sec:general} is dedicated to showing the existence of a Gaussian process with mean and covariance defined as above but in a more general setting.

\section{General Results}\label{sec:general}

Recall the form of a Gaussian process $\mathbb{X} = \{X_s: \, s \in T\}$ whose value is fixed at several points $\{t_1, ..., t_n\}$, and whose mean and covariance are given by equations \eqref{eq:finiteMean} and \eqref{eq:finiteKernel}, respectively. Section \ref{sec:ortho} provides a construction for such a process using the theory of RKHS. In this section we apply the same framework when $T_0$ is an arbitrary compact subset of the input domain. Our approach is to first show that such a Gaussian process indeed exists and can be described using only the information on $T_0$. We then show that this process can be arrived at by taking the limits of \eqref{eq:finiteMean} and \eqref{eq:finiteKernel} when the collection of points $\{t_1,...,t_n\}$ approaches a dense subset of $T_0$.

\subsection{Construction} \label{sec: constr}

Let $T\subset \mathbb{R}^d$ and $T_0 \subset T$ {be compact sets}, $k$ a continuous {and universal} covariance kernel on $T$ \cite{micchelli}, and $g$ an element of $\mathcal{H}(T)$. First observe that any Gaussian process which is fixed on $T_0$ must have a covariance function $k_0$ satisfying $k_0(s,t)=0$, if either $s \in T_0$ or $t \in T_0$; i.e., the desired covariance kernel must vanish on $T_0 \times T$. 

Let $\mathcal{H}_0 = \{f \in \mathcal{H}(T):\, f|_{T_0} \equiv 0\}$ as in Subsection \ref{sec:ortho}. Since 
{$\mathcal{H}_0$ is a closed subspace of $\mathcal{H}(T)$}, there exists an orthogonal projection $P: \mathcal{H}(T) \to \mathcal{H}_0$ and hence a kernel $k_0 = Pk$ that generates $\mathcal{H}_0$.  {We require that} the mean of the conditional distribution equals $g$ on $T_0$.  Thus, define $[g] = \{ f \in \mathcal{H}(T): \mbox{$f({T_0})=g({T_0})$}\}$.  For $f\in \mathcal{H}(T)$, let $f=f_0+f_{\perp}$ be the unique decomposition of $f$ with $f_0 \in \mathcal{H}_0$ and $f_{\perp} \in \mathcal{H}_0^\perp$.  Note that $f\in [g]$ if and only if $f-g \in \mathcal{H}_0$, which in turn is true if and only if $f_\perp = g_\perp$.  In other words, $[g] = \{ f\in \mathcal{H}(T): f_\perp = g_\perp \}$, and our requirement on the conditional mean function is that it belongs to $[g]$.

The Kolmogorov Existence Theorem permits the existence of a Gaussian process given a mean $\mu$ and kernel function $k$ provided that the $k$ is symmetric and positive semi-definite \cite[Theorem 5.16]{kall}. As a corollary, we have the following result.

\begin{theorem}
Given a continuous, symmetric, positive semi-definite covariance function, $k$, and $\mu \in \mathcal{H}(T)$, there exists a Gaussian process $\mathbb{X}=\{X_t: \, t \in T\}$ with mean $\mu_0 = P \mu + g_\perp$, covariance {kernel} $Pk$, {and such that} $X_t = g_\perp(t)$ (a.s.) for each $t \in T_0$.
\end{theorem}

It remains to see how one might compute $P\mu$ for arbitrary $\mu \in \mathcal{H}(T)$. {Similar to the technique used in Section \ref{sec:rkhs_intro}, consider the behavior of elements of $\mathcal{H}(T)$ restricted to $T_0$. We will show that there is an equivalence between $\mathcal{H}(T_0)$ and $\mathcal{H}_0^\perp$. The following proposition is important in that it provides a useful characterization of $\mathcal{H}_0^\perp$.}
\begin{proposition} \label{prop: perp}
    $\mathcal{H}_0^\perp = \overline{\mathrm{Span}}(\{k_s: \, s \in T_0\}).$

    \begin{proof}
    
    ($\supseteq$) Note that for any $s \in T_0$, and any $f \in \mathcal{H}_0$, $\langle k_s, f \rangle_{\mathcal{H}(T)} = f(s) = 0$. Thus, 
    $\{k_s: \, s \in T_0\} \subset \mathcal{H}_0 ^ \perp$, which implies $\mathcal{H}_0^\perp \supset \overline{\mathrm{Span}}(\{k_s: \, s \in T_0\})$. 

    ($\subseteq$) It suffices to show that $\mathcal{H}_0 \supset \overline{\mathrm{Span}}(\{k_s; \, s \in T_0\})^\perp.$ Let $f \in \overline{\mathrm{Span}}(\{k_s: \, s \in T_0\})^\perp$. For any $s \in T_0$, $\langle f, k_s \rangle_{\mathcal{H}(T)} = 0$, implying that $f(s) = 0$ and, hence, $f|_{T_0} \equiv 0$.
    \end{proof}
\end{proposition}

{As with $\mathcal{H}(T)$, we can obtain an RKHS of functions on $T_0$ via closing the span of the restricted functionals $k_s|_{T_0}, ~s \in T_0$; i.e.,} $\mathcal{H}(T_0) = \overline{\mbox{Span}}(\{k_s|_{T_0}; \, s \in T_0\})$. Hence, the equivalence between $\mathcal{H}(T_0)$ and $\mathcal{H}_0^\perp$ {can be established via unique extension of} each element of $\mathcal{H}(T_0)$ to all of $T$. 

\begin{theorem}\label{thm:iso_thm}
There exists an isometric isomorphism between $\mathcal{H}_0^\perp$ and $\mathcal{H}(T_0)$.

\begin{proof}
    See Appendix \ref{app:thm_3.3}.
\end{proof}
\end{theorem}

{For ease of notation}, { for $ f \in \mathcal{H}(T)$, set $\hat{f} = (\Tilde{\psi} \circ Q) f$ where $Q : \mathcal{H}(T) \rightarrow \mathcal{H}_0^\perp$ projects into $\mathcal{H}_0^\perp$ and  $\Tilde{\psi}: \mathcal{H}_0^\perp \to \mathcal{H}(T_0)$ is the isometric isomorphism defined in the proof of Theorem \ref{thm:iso_thm} that maps a function $\tilde{f} \mapsto \tilde{f}|_{T_0}$. Then, again assuming as we did in subsection \ref{sec:ortho} that $g \in \mathcal{H}_0^\perp$, Theorem \ref{thm:iso_thm} and equation \eqref{eq: mean_projection} yield }
{
\begin{align}
    \mu_0(s) &= \mu(s) ~{+}~ {Q}(g - \mu)(s) \notag \\
        &= \mu(s) ~{+}~ \langle {Q}k_s, {Q}(g - \mu) \rangle_{\mathcal{H}(T)} \notag \\
        &= \mu(s) ~{+}~ \langle\hat{k}_s, \hat{g} - \hat{\mu} \rangle_{\mathcal{H}(T_0)}, \label{eqn:condMeanfun}
\end{align}
where the second line follows from $Q$ being self-adjoint and idempotent, and the last line follows from the fact that $\tilde{\psi}$ is 
an isometry.} Similarly, by \eqref{eq: cov_projection},
\begin{align}
    k_0(s_1, s_2) &= k(s_1, s_2) - {Q}k(s_1, s_2) \notag \\
        &= k(s_1, s_2) - \langle {\hat{k}_{s_1}}, {\hat{k}_{s_2}} \rangle_{\mathcal{H}(T_0)}. \label{eqn:condCovfun}
\end{align}
{$\mu_0(\cdot)$ and $k_0(\cdot, \cdot)$ are the mean function and the kernel function that define our proposed projected kernel Gaussian process (pkGP). In other words, starting with a typical Gaussian process $\mathcal{GP}(\mu(\cdot), k(\cdot, \cdot))$, pkGP is the Gaussian process $\mathcal{GP}(\mu_0(\cdot), k_0(\cdot, \cdot))$ where $\mu_0$ and $k_0$ are defined by \eqref{eqn:condMeanfun} and \eqref{eqn:condCovfun}, respectively.} {Approximating the inner products $\langle\hat{k}_t, \hat{g} - \hat{\mu} \rangle_{\mathcal{H}(T_0)}$ and $\langle \hat{k}_{s_1}, {\hat{k}_{s_2}} \rangle_{\mathcal{H}(T_0)}$ is discussed in Section \ref{sec:rrAndKrig}. Here we remark that doing so only requires knowing the kernel sections and the functions $\mu$ and $g$ on $T_0$ so that we can take the needed $L^2$ inner products on $T_0$. These are of course assumed known since we are conditioning the GP on $g(T_0)$.}

\paragraph{Remarks} {While we assume that $T_0$ is compact throughout this paper, this assumption is not necessary for Theorem \ref{thm:iso_thm} to hold. We also remark that results nearly identical to Proposition \ref{prop: perp} and Theorem \ref{thm:iso_thm}  were derived by \cite[Lemma B.2]{chen2025} simultaneously and independently of our work, but in the context of quasi Monte Carlo methods for approximating integrals. The claims of \cite{chen2025} and our paper are both reminiscent of those appearing in \cite{GauthierBay12}. The different contexts in which the connection between $\mathcal{H}_0$ and $\overline{\mathrm{Span}}(\{k_s; \, s \in T_0\})$ appears underscores its broad utility.}

\subsection{Connection to Finite Case} \label{sec:limits}
The purpose of this subsection is to connect the finite-dimensional case in Subsection \ref{sec:ortho} to the general case in Subsection \ref{sec: constr} by showing that the same result can be obtained by taking limits of countable sets. An interpretation of this is that if one selects enough points on $T_0$ as pseudo-training data for the Gaussian process, the resulting GP conditioned on these points serves as a justifiable approximation to pkGP with mean and covariance given by \eqref{eqn:condMeanfun} and \eqref{eqn:condCovfun}. Assuming one selects an appropriate subset of $T_0$, this can be computed to arbitrary precision. {In Section \ref{sec:simulations} we empirically demonstrate this claim.}

It is important that one can construct a countable dense subset of $T_0$. (This is the case if $T_0 \subset \mathbb{R}^d$.) By Proposition \ref{prop: perp} {and Theorem \ref{thm:iso_thm}}, we need not worry about considerations of $T_0$ as a subset of $T$, but rather can simply consider $T_0$ as its own space. Hence, we assume any function defined in this subsection is done so on $T_0$. Let $D= \{t_n\}$ be a countably dense subset of $T_0$, and consider $\mathcal{K}_D:=\overline{\mbox{Span}}(\{{\hat{k}_t}: \, t \in D\})$. Since $D$ is dense, for arbitrary $s \in T_0$, there exists a subsequence $\{t_{n_j}\} \subset D$ so that ${\hat{k}_s} = \lim_{j \to \infty} {\hat{k}_{t_j}}$, {by continuity of the kernel sections (Lemma \ref{lemma:kHolder})}. Therefore,
$$
{\mbox{Span}}\{{\hat{k}_s}; \, s \in T_0\} \subset \mathcal{K}_D \subset \mathcal{H}(T_0).
$$
Taking closure, we see that $\mathcal{K}_D=\mathcal{H}(T_0)$. 

A consequence of the preceding is that, for a given ${\hat{f}} \in \mathcal{H}(T_0)$ and for $\epsilon > 0 $, there exists an $N_0$ so that any interpolating approximation $f_N$ of {$\hat{f}$} {obtained from a finite subset} $\{k_{t_n}\}_{n=1}^N$ satisfies
$$
||f_N- {\hat{f}}||_{\mathcal{H}(T_0)} < \epsilon, \, \mbox{for } N \ge N_0.
$$
By defining {$Q_N : \mathcal{H}(T_0) \rightarrow {\overline{\mbox{Span}}}(\{\hat{k}_{t_n}\}_{n=1}^N)$ so that $Q_N{\hat{f}} = f_N$}, this is equivalent to saying that $Q_N$ converges strongly to the identity operator. Strong operator convergence implies that for ${\hat{f}, \hat{g}} \in \mathcal{H}(T_0)$, $\langle Q_N{\hat{f}, \hat{g}} \rangle_{\mathcal{H}(T_0)} \to \langle Q{\hat{f}, \hat{g}} \rangle_{\mathcal{H}(T_0)}$. {Convergence of RKHS interpolators as mentioned here have been studied extensively in the literature (e.g., \cite{Traub88, chen2025})}. In our case, however, we have that the finite dimensional mean and covariance from Subsection \ref{sec:ortho} converges to the infinite dimensional mean and covariance from Subsection \ref{sec: constr} {as the pseudo-observation points in $T_0$ approach a dense subset}.

\subsection{Weak Convergence of the Stochastic Process}\label{sec:prob}

The previous section finds, under certain conditions, closed-form expressions for the mean and covariance of the proposed pkGP. In addition, the preceding subsection provides a means of {reliably approximating} the mean and covariance via selecting a representative finite subset of $T_0$. The aim of this subsection is to show that this limiting approximation of the mean and covariance can be leveraged to establish weak convergence of the process itself. 

Define $\mu_0^N$ and $k_0^N$ to be the mean and covariance functions, respectively, resulting from conditioning on $\{t_1,\hdots,t_N\}, ~N < \infty$. (See \eqref{eq:finiteMean} and \eqref{eq:finiteKernel}.) Let  $\{\mathbb{X}^N\}_{N=1}^\infty$ be {the} sequence of GPs such that $\mathbb{X}^N$ has mean $\mu_0^N$ and covariance $k^N_0$, and define $\mathbb{X}$ to be the limiting pkGP with mean and covariance $\mu_0$ and $k_0$ {defined in equations \eqref{eqn:condMeanfun} and \eqref{eqn:condCovfun}, respectively}.

Subsection \ref{sec:limits} establishes the convergence of $\mu^N$ and $k^N$ to $\mu_0$ and $k_0$, respectively. Thus, the convergence of any finite dimensional distribution of $\mathbb{X}^N$ to that of $\mathbb{X}$ is obtained. To show weak convergence {via Prohorov's Theorem}, it remains to show that {the sequence} $\{\mathbb{X}^N\}_{N=1}^\infty$ is tight {\cite[Theorem 7.1]{billConvergence}}. We first provide conditions under which one may find a version of $\mathbb{X}$ which is continuous, as this is a desirable property in practice. The proofs of {Lemmas \ref{lem:continuous} and \ref{lem:tight}} may be found in Appendix A.

\begin{lemma}\label{lem:continuous}
Suppose that $\mathbb{X}$ is a Gaussian process with mean $\mu$ and covariance kernel $k$. If $\mu$ is continuous and $k$ is $\gamma-$H\"older continuous on $\mathbb{R}^d \times \mathbb{R}^d$, then there is a version of $\mathbb{X}$ which almost surely continuous.
\begin{proof}
    {This is a direct consequence of \cite[Theorem 2.14]{Nummi2024}.} (See also \cite{dacosta2024}.) For completeness, though, we provide a proof in Appendix \ref{app:lem_3.4}.
\end{proof}
\end{lemma}

It is indeed the case that $\{\mathbb{X}^N\}_{N=1}^\infty$ is tight if the conditions for the Kolmogorov-Chentsov theorem stated above are met uniformly on $N$. This is established in the following {lemma}.

\begin{lemma}\label{lem:tight}
If the covariance kernel $k$ is $\gamma-$H\"older continuous, $k$ is universal on $T_0$ and $g|_{T_0},\mu|_{T_0} \in \mathcal{H}(T_0)$, then $\{\mathbb{X}^N\}_{N=1}^\infty$ is tight in $(C(T),||\cdot||_\infty)$.
\begin{proof}
    See appendix \ref{app:thm_3.5}.
\end{proof}
\end{lemma}
{Lemmas \ref{lem:continuous} and \ref{lem:tight} immediately yield weak convergence \cite[Theorem 14.3]{kall}:}
\begin{theorem}\label{thm:weakConv}
    {Under the conditions of Lemmas \ref{lem:continuous} and \ref{lem:tight}, $\mathbb{X}^N \stackrel{w}{\rightarrow} \mathbb{X}$ in $C(T)$}.
\end{theorem} 
{\paragraph{Remarks} Here it is worth mentioning that we are using interchangeably the notions of tightness of the Gaussian measures on $C(T)$ and tightness of the sequence of random processes $\{\mathbb{X}^N\}_{N=1}^\infty$. The latter follows the view of \cite{kall}, who uses the notion in defining and proving ``weak convergence" of the sequence of processes. Sharing this view, this is what we mean by ``weak convergence" in this section. Further, we remark that while the GPs having continuous versions is a desirable property, it is not continuity that directly leads to tightness. It is the antecedent conditions that are common to both the continuity and convergence results. The interested reader may find results similar to Lemma \ref{lem:continuous}, namely when GPs lie in so-called reproducing kernel Banach space, in \cite[Appendx B]{pfortner2024}.}

\subsection{Practical Considerations: Including a Nugget}\label{subsec:nugget}

Here we briefly consider a Gaussian process modification that is often used in practice: the addition of a nugget. We discuss how this relates to our proposed approach.

A common use of GP models is for emulating deterministic computer output from a complex and computationally expensive model \cite{SantnerEtAl03}. In other words, we are taking computer model input/output $\{(t_i, y_{t_i})\}$ and training a GP to interpolate these points in some bounded subset $T \subset \mathbb{R}^d$. GP emulators for computer models commonly employ the squared exponential covariance kernel, defined by $k(s,t) = \exp\{-\sum_{k = 1}^d \ell_k^{-1}|s_k - t_k|^2\}$, where $s, t \in \mathbb{R}^d$, and $\ell_1, \ldots, \ell_d > 0$. This covariance function produces very smooth sample paths at the cost of a poorly conditioned covariance matrix. It is therefore commonplace when using this kernel to employ a regularization component to bound the condition number, referred to as a ``nugget" \cite{ranjan11}, thereby improving the stability of matrix computations. In this case, a covariance matrix of the form $k(\mathbf{s}, \mathbf{t}) = (k(s_i, t_j))_{i, j}$ instead becomes $k(\mathbf{t}, \mathbf{t}) + \delta \mathbf{I}$, where $\mathbf{I}$ is the identity matrix, and $\delta$ is a small number that can be tuned through a variety of means \cite{ranjan11}. This results in a process which is, strictly speaking, no longer continuous. In practice, though, the sample paths are for most purposes nearly identical to those arising from the original process without a nugget. In addition, there is often practical justification for this nugget to be added; e.g., to represent measurement error or to improve predictive ability \cite{GramacyLee12}.

Similar computational considerations can be made with our proposed approach. In keeping with our practice of avoiding direct matrix operations, let us consider the linear operator mapping $L^2(T_0)$ to itself defined by $\Tilde{K} = K + \delta I$, where $K$ is defined in \eqref{eq:integral_operator}, and $I$ is the identity operator. Recalling a form of the RKHS inner product provided in Subsection \ref{sec:rkhs_intro}, we have
$$
\langle {\hat{f},\hat{g}} \rangle_{\mathcal{H}(T_0)} = \langle K^{-1/2}{\hat{f}}, K^{-1/2}{\hat{g}} \rangle_{T_0},
$$
where $\langle \cdot, \cdot \rangle_{T_0}$ denotes the $L^2$ inner product on $T_0$. Using the same notation as in Section \ref{sec:prelim}, the eigenvalues and eigenvectors of $\Tilde{K}$ are $\{\lambda_n+ \delta\}$ and $\{e_n\}$, and so one may represent $\Tilde{K}$ as $
\tilde{K}(\cdot)= \sum_{n=1}^\infty (\lambda_n+\delta) \langle \,\cdot \,,e_n\rangle_{T_0} e_n.
$
The eigenvalues of $\Tilde{K}$ are bounded below by $\delta$, implying that $\Tilde{K}$ has a bounded inverse operator $\Tilde{K}^{-1}$. Therefore $\Tilde{K}^{-1/2}$ can be represented by $\tilde{K}^{-1/2}(\cdot)= \sum_{n=1}^\infty (\lambda_n+\delta)^{-1/2}\langle \,\cdot \,,e_n\rangle_{T_0} e_n$.
Replacing $K^{-1/2}$ in \eqref{eq:rootInnerProduct} with $\Tilde{K}^{-1/2}$, we obtain an approximation for the RKHS inner product for ${\hat{f}_1}, {\hat{f}_2} \in \mathcal{H}(T_0)$ as 
\begin{equation}\label{eqn:nuggIP}
\langle {\hat{f}_1}, {\hat{f}_2} \rangle_{\Tilde{K}} = \langle \Tilde{K}^{-1/2}{\hat{f}_1},\Tilde{K}^{-1/2} {\hat{f}_2} \rangle_{T_0} = \sum_{n=1}^\infty \frac{\langle {\hat{f}_1},e_n\rangle_{T_0} \langle {\hat{f}_2},e_n \rangle_{T_0}}{\lambda_n + \delta}.
\end{equation}
It follows that the pkGP with posterior mean $\Tilde{\mu}_0$ and posterior covariance $\tilde{k}_0$ may be obtained by replacing $\langle \cdot, \cdot \rangle_{\mathcal{H}(T_0)}$ with $\langle \cdot, \cdot \rangle_{\Tilde{K}}$ in \eqref{eq: cov_projection} and \eqref{eqn:condCovfun}; i.e., 
\begin{align}\label{eqn:nuggFnArray}
    \Tilde{\mu}_0(s_1) & = \mu(s_1) + \langle {\hat{k}_{s_1}}, {\hat{g}}-{\hat{\mu}}\rangle_{\Tilde{K}}, \notag \\
    \Tilde{k}_0(s_1,s_2) & = k(s_1,s_2) - \langle {\hat{k}_{s_1}}, {\hat{k}_{s_2}} \rangle_{\Tilde{K}}.
\end{align}
Again, this process will no longer have continuous sample paths. Assuming $\delta$ is small enough, though, this is not an obstacle in practice.

\section{Computing RKHS Inner Products}\label{sec:rrAndKrig}
The previous sections show that one may construct a Gaussian process which has zero variation on an arbitrary subset $T_0$ of the domain, and define its mean and covariance functions in terms of an RKHS inner product. In practice, however, the RKHS inner product cannot be computed exactly. Here we discuss a technique for {computing RKHS inner products via Rayleigh-Ritz approximation of the associated kernel spectra, and compare it to the more direct approach of SVD-based kriging interpolation}. {To illustrate, we consider both a function that is contained in the associated RKHS, and one that is not.}

\subsection{Computation of RKHS Inner Product}\label{sec:computation}
Recall that the RKHS norm is given in terms of the spectral decomposition $\{(\lambda_n, e_n) \}$ of the integral operator {$K$}, which must be computed numerically. The inner product $\langle \cdot, \cdot \rangle_{\mathcal{H}(T_0)}$ {then may be} approximated via the bilinear form $a_N ( \cdot, \cdot )$, given by
$$
a_N({\hat{f}}, {\hat{g}} ) = \sum_{n=1}^N \frac{\langle {\hat{f}}, e_n \rangle_{T_0} \langle {\hat{g}}, e_n \rangle_{T_0} }{\lambda_n}.
$$

While the form of $a_N (\cdot, \cdot )$ does not permit a convergence independent of the choice of ${\hat{f}, \hat{g}} \in \mathcal{H}(T_0)$, {uniform} convergence can be established for the family of functions $\mathcal{K}:=\{k_t: t \in T_0\}$.
\begin{proposition}
The {sequence} of bilinear forms $\{a_N\}$ converges uniformly to $\langle \cdot, \cdot \rangle_{\mathcal{H}(T_0)}$ on $\mathcal{K} \times \mathcal{K}$ {as $N \rightarrow \infty$}.
\begin{proof}
Define $F_N, F: T_0 \times T_0 \to \mathbb{R}$ by $F_N(s,t) = a_N({\hat{k}_s}, {\hat{k}_t})$ and $F(s,t) = \langle {\hat{k}_s}, {\hat{k}_t} \rangle_{\mathcal{H}(T_0)}$. It is clear that $F_N \to F$ pointwise. {Hence, by the Arzel\'{a}-Ascoli Theorem,} it suffices to show that $\{F_N\}$ is equicontinuous. Defining $Q_N$ to be the projection from $\mathcal{H}(T_0)$ to $\mbox{Span}( \{e_n\}_{n=1}^N)$, we have that 
$$
F_N(s,t) = \langle Q_N {\hat{k}_s}, Q_N {\hat{k}_t} \rangle_{\mathcal{H}(T_0)},
$$ 
and so equicontinuity follows directly from the fact that $F$ is H\"older continuous and $\{Q_N\}$ is uniformly bounded by the identity operator.
\end{proof}
\end{proposition}
{The RKHS inner product is defined and evaluated via the eigensystem of the integral operator, $\{(\lambda_n, e_n)\}_{n=1}^N$, which can be difficult to compute directly. \cite{oya2009} propose to use the Rayleigh-Ritz (RR) method to approximate the eigenvectors and associated eigenvalues, whence the RKHS inner product can be approximated to arbitrary precision.} 

The RR approach proceeds by first selecting $m$ linearly independent functions $\{\xi_j\}_{j=1}^m$ on $T_0$, whence the eigenfunctions are approximated as $\tilde{e}_i = \sum_{j=1}^m b_{ij}\xi_j, ~i= 1, \ldots, {m}$. The coefficients $\boldsymbol{b}_i = (b_{i1}, \ldots, b_{im})^\top$ and approximate eigenvalues $\{\tilde{\lambda}_i\}_{i=1}^m$ are obtained by solving the generalized eigenvalue problem,
\begin{equation}\label{eqn:geneval}
    \boldsymbol{B}\boldsymbol{b}_i = \tilde{\lambda}_i\boldsymbol{C}\boldsymbol{b}_i, ~~i= 1, \ldots, m,
\end{equation}
where $(\boldsymbol{B})_{ij} = \langle K\xi_i, \xi_j \rangle$ and $(\boldsymbol{C})_{ij} = \langle \xi_i, \xi_j \rangle$. These $L^2$ inner products can be evaluated via quadrature {over $T_0$}. {(In this paper we use Gauss-Legendre quadrature.)} Observe that when the set $\{\xi_j\}_{j=1}^m$ is chosen to be orthonormal (e.g., orthogonal polynomials), $\boldsymbol{C} = \boldsymbol{I}$ and \eqref{eqn:geneval} becomes an ordinary eigenvalue problem.

With $(\tilde{\lambda}_i, \tilde{e}_i), ~i= 1, \ldots, m$ in hand, {the inner products of interest are approximated with 
\begin{align}
    \langle \hat{f}, \hat{g} \rangle_{\tilde{H}_n(T_0)} &:= \sum_{i=1}^n \tilde{\lambda}_i^{-1} \langle \hat{f}, \tilde{e}_i \rangle_{T_0} \langle \hat{g}, \tilde{e}_i \rangle_{T_0}, ~~\hat{f}, \hat{g} \in \mathcal{H}(T_0)
    \label{eqn:oyaInner}
\end{align}
where $n \leq m$.} The validity of these approximations is established via the following theorem:
\begin{theorem}
    {For $\langle \cdot, \cdot \rangle_{\tilde{H}_n(T_0)}$ as defined above, and for $\hat{f}, \hat{g} \in \mathcal{H}(T_0)$, 
    $$
        |\langle \hat{f}, \hat{g} \rangle_{\tilde{H}_n(T_0)} - \langle \hat{f}, \hat{g} \rangle_{\mathcal{H}(T_0)}| \rightarrow 0
    $$
}
{as $n, m \rightarrow \infty$.} 
\end{theorem}
\begin{proof}
    \cite[Appendix A]{oya2009}.
\end{proof}
The choice of basis functions may depend on the specific application. For instance, polynomials for smoothly-varying processes, or wavelets for non-smooth covariance functions.

\subsection{Numerically Verifying the Reproducing Property}\label{subsec:verify}
{{For an input domain $T_0$} and a given function $f \in \mathcal{H}(T_0)$, the reproducing property 
\begin{equation}\label{eqn:interp}
    f(t) = \langle f, k_t \rangle_{\mathcal{H}({T_0})}
\end{equation} leads to the RKHS interpolator. {(In this subsection, for simplicity we have dropped the hat notation for functions in $\mathcal{H}(T_0)$ so that the notation can be used for estimators, following convention.) }{In the statistics literature, the use of Gaussian processes leads to the kriging interpolator, 
\begin{equation}\label{eqn:kriging}
    \hat{f}_{krig}(\cdot) = \boldsymbol{k}_{\boldsymbol{x}}(\cdot)\boldsymbol{K}^{-1}_{\boldsymbol{x}, \boldsymbol{x}}\boldsymbol{f}(\boldsymbol{x}).
\end{equation}}
{The kriging predictor has its roots in spatial statistics and is known to be the optimal predictor among linear unbiased estimators. Despite the differences in motivation and interpretation, kriging is essentially equivalent to the RKHS interpolant in \eqref{eqn:interp}, which itself is known to be optimal in both the RKHS norm and among linear estimators \cite{Wend05}. This can be seen by observing that $\hat{f}_{krig}(t) = (\boldsymbol{K}_{\boldsymbol{x}, \boldsymbol{x}}^{-1/2}\boldsymbol{k}_{\boldsymbol{x}}^\top(t,\boldsymbol{x}))^\top(\boldsymbol{K}_{\boldsymbol{x}, \boldsymbol{x}}^{-1/2}\boldsymbol{f}(\boldsymbol{x})) \approx \langle K^{-1/2}k_t, K^{-1/2}f\rangle_{T_0} = \langle f, k_t\rangle_{\mathcal{H}(T_0)}$.}

{In practice, kernel functions such as the popular squared exponential kernel yield poorly conditioned covariance matrices, meaning the inverse of $\boldsymbol{K}_{\boldsymbol{x},\boldsymbol{x}}$ in \eqref{eqn:kriging} is computationally unstable. As discussed in subsection \ref{subsec:nugget}, this is typically handled by adding a small nugget to $\boldsymbol{K}_{\boldsymbol{x},\boldsymbol{x}}$ before inversion or singular value decomposition. This approximation of course comes at a cost in terms of how well the eigenvalues of $K$ are estimated. On the other hand, the Rayleigh-Ritz algorithm is such that $\boldsymbol{K}_{\boldsymbol{x},\boldsymbol{x}}$ does not need to be computed directly. It can thus often avoid the use of a nugget at all, even in otherwise poorly conditioned cases like the squared exponential. In this section, we briefly illustrate the consequences of using a nugget by showing improvement in approximation error that can be achieved when one avoids the nugget entirely via RR. While it is widely understood that nuggets change an approximation, we include a direct comparison here since, to the best of our knowledge, such comparisons are rarely done and the Rayleigh-Ritz method has not appeared much in the statistics literature \cite{ebner24}.}

{We take the {domain of interest to be $T_0 = [-1, 1]$} and the kernel to be the squared exponential (Gaussian) kernel, $k(x, x^\prime) = \exp\{-|x - x^\prime|^2\}$. We can create a target function in the associated RKHS {$\mathcal{H}(T_0)$} with
\begin{equation}\label{eqn:frkhs}
    f_{_{RKHS}}(\cdot) = \sum_{i=1}^4 \alpha_i k(\cdot, x_i),
\end{equation} 
where $\alpha_1, \ldots, \alpha_4$ are drawn independently from $\text{Unif}(-1, 1)$ and the $x_1, \ldots, x_4$ are regularly spaced between -1 and 1. To construct a function that is {\em not} in the RKHS, we take $f_{_{Lagrange}}(\cdot)$ to be the Lagrange polynomial that interpolates $f_{_{RKHS}}$ given a set of interpolation points. Since the Langrange function is a polynomial, it cannot be a member of the RKHS \cite{ha}. Instances of two such functions are plotted in Figure \ref{fig:rkhsInterpEx}, in which they can be seen to be very similar but not equal. To approximate the eigenfunctions in the RR algorithm, we use Legendre orthogonal polynomials; i.e., $\tilde
{e}_i(\cdot)= \sum_{j=1}^M a_{ij}\xi_j(\cdot)$ and \eqref{eqn:geneval} becomes an ordinary eigenvalue problem. For a fair comparison using the same amount of information, we use $M$ evenly spaced observations, $\boldsymbol{x}^{(M)} = (x_1, \ldots, x_M)^\top$, as training points for the kriging predictor \eqref{eqn:kriging}. {Our implementation of SVD toward inverting $\boldsymbol{K}^{-1}_{\boldsymbol{x}, \boldsymbol{x}}$ in \eqref{eqn:kriging} is via the \texttt{eigen} function in base \texttt{R}, which uses the LAPACK package \cite{Lapack99}. We add a nugget of $10^{-6}$ prior to the SVD \eqref{eqn:kriging}. We do {\em not} use a stabilizing nugget for the Rayleigh-Ritz (RR) algorithm.}
\begin{figure}[tb]
    \centering
    \includegraphics[scale= 0.5]{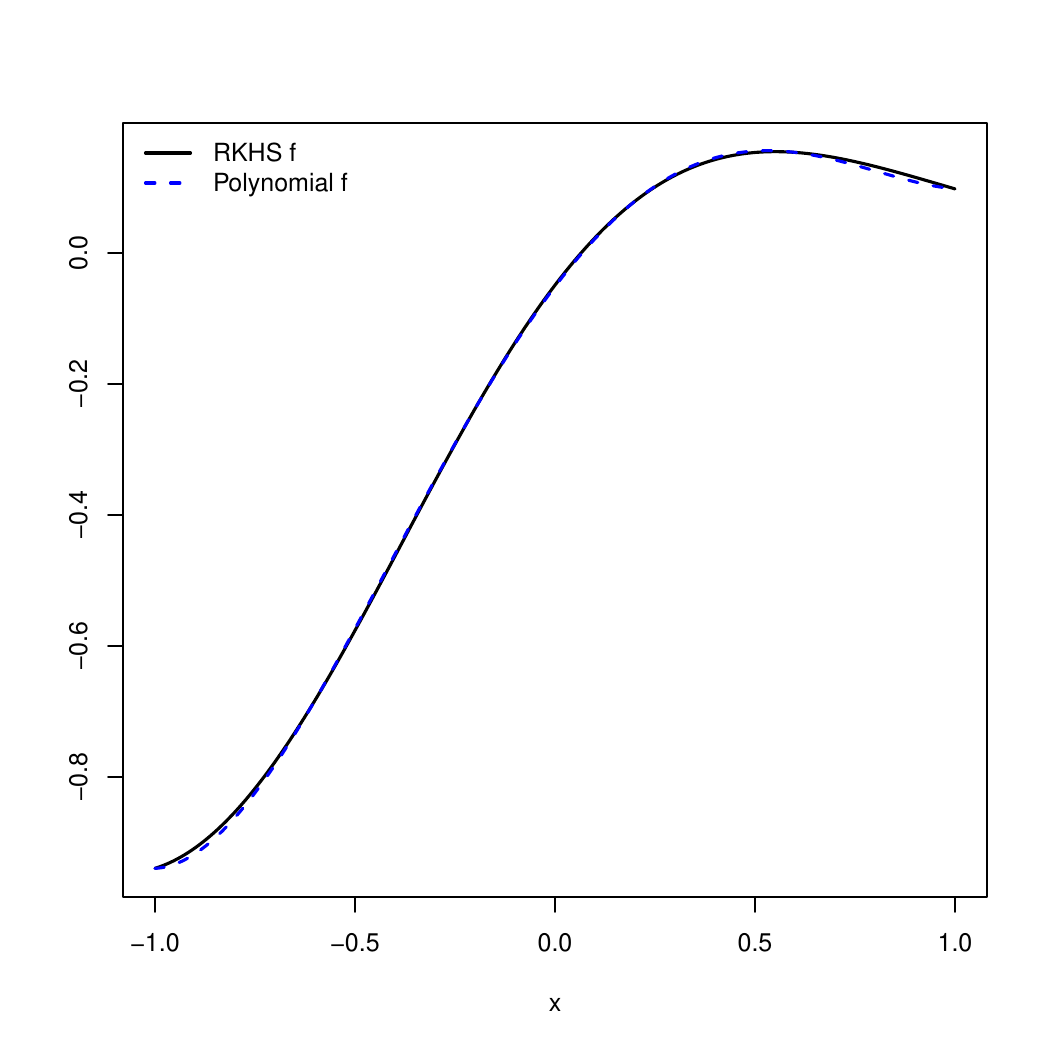}
    \caption{{\em Two example functions to be reproduced via interpolation. One is an element of the RKHS, and the polynomial is not.}}
    \label{fig:rkhsInterpEx}
\end{figure}

{Figure \ref{fig:InterpErrorEx} displays the relative errors with respect to the supremum norm, defined as $\|\hat{f} - f \|_{\infty}/\|f\|_\infty$. In both plots, we vary the RR eigenfunction approximation order $M$ from 6 to 55. The curves are calculated over 100 random functions simulated according to \eqref{eqn:frkhs}. We observe that when reproducing a function that is a member of the RKHS, the RR approximation vastly outperforms the {SVD} calculation, regardless of the value of $M$. Further, even when the target function is not an element of the RKHS, which is most likely the case in practice, RR-based estimation results in a reconstruction considerably better than {SVD-based kriging}. The difference becomes more pronounced as $M$ increases, though the RR approximation error at high values of $M$ pulls RR back toward SVD-based kriging. {Our illustrative example thus suggests that for the commonly used Gaussian kernel, approximating the RKHS inner product via RR spectral decomposition might be the preferred approach.}
\begin{figure}[tb]
    \centering
    \includegraphics[scale= 0.4]{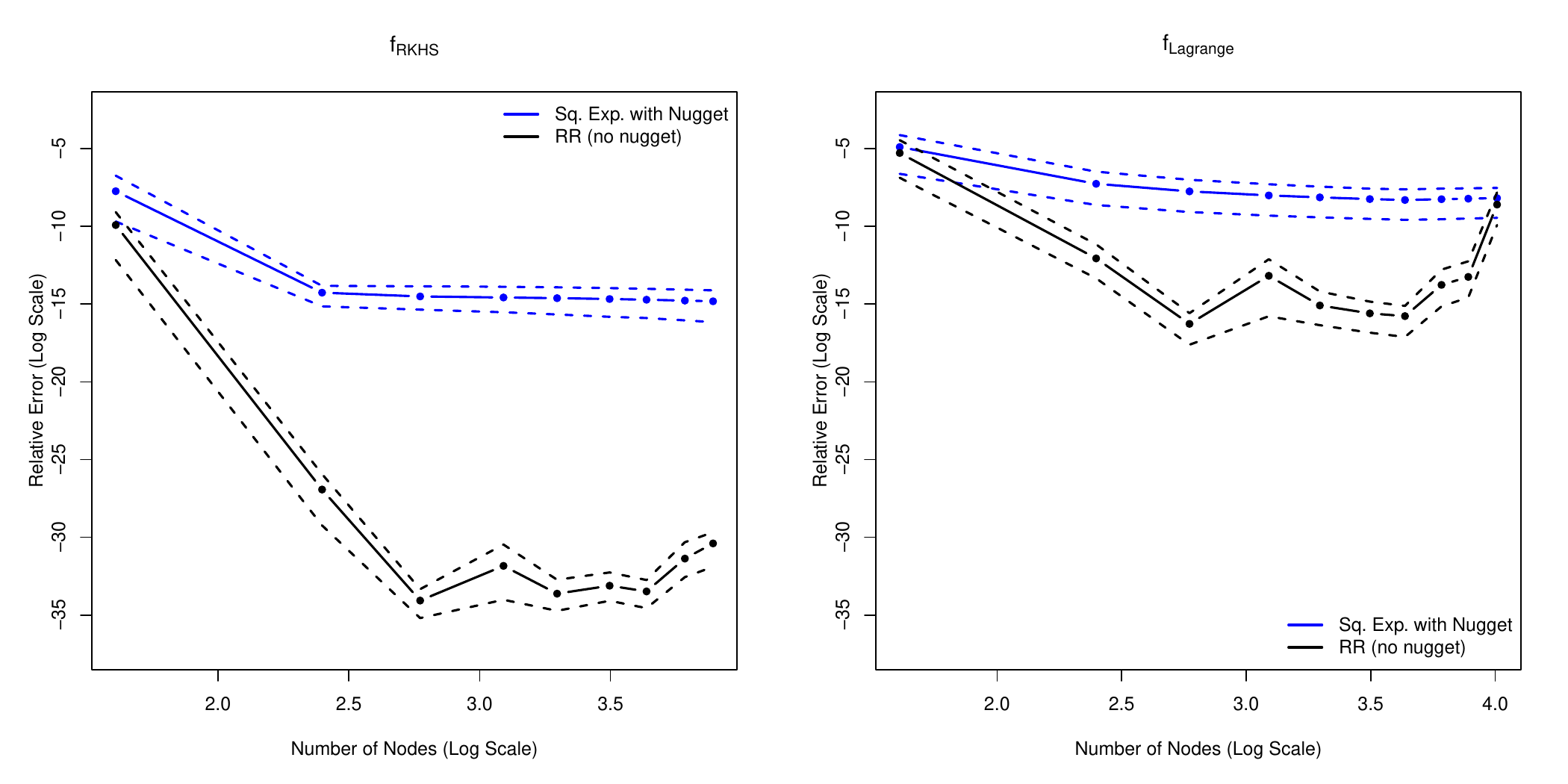}
    \caption{{\em Median reproduction errors of the {SVD-based kriging interpolator (squared exponential with a nugget) \eqref{eqn:kriging} vs. the RR-based approximation of \eqref{eqn:interp} (no nugget)} over 100 replications. The left panel is for a function that is an element of the RKHS. The right panel is for the polynomial function that is not. The dashed bands about each curve denote the empirical pointwise 95\% confidence intervals.}}
    \label{fig:InterpErrorEx}
\end{figure}

\section{Numerical Experiments}\label{sec:simulations}
Here we consider simulated examples to illustrate our proposed approach in two different scenarios. The first case is that of known boundary conditions; e.g., Dirichlet conditions for a system of partial differential equations. In this situation we compare our approach to the na\"{i}ve approach of adding a finite number of (known) function values along on the boundary as pseudo-training data, as this is what a practitioner might do since these pseudo-data are available at no additional cost. The second example we consider is one in which the function is not known along the boundary of the domain, but rather along a subset tracing a diagonal of the domain. 

{For clarity prior to presenting the results, we summarize the procedure for constructing the projected mean and covariance functions in Algorithm \ref{ALG:pkGPAlg}. These projected functions determine the pkGP which can be used as the prior GP for training after collecting training points. Throughout this section, we assume the base Gaussian process model (prior to projection) has mean zero, $\mu(t) \equiv 0$, as mean zero GP's have been shown to often suffice for surrogate modeling tasks \cite{Bayarri07}.}
\begin{algorithm}[tb]
{
    \caption{{Sketch of the algorithm for constructing the projected kernel Gaussian process (pkGP). Note that the nugget $\delta$ often may be taken to be zero.}}
    \label{ALG:pkGPAlg}
    \begin{algorithmic}[1]
    \REQUIRE GP mean function $\mu(\cdot)$ and covariance function $k(\cdot, \cdot)$ with regularizing nugget $\delta$, conditioning domain subset $T_0$, known function $g(T_0)$
    \STATE Approximate the spectrum of the integral operator $K$ using Rayleigh-Ritz; i.e., choose basis functions $\{\xi_j\}_{j=1}^m$ and a quadrature rule to compute the matrices $\boldsymbol{B}, \boldsymbol{C}$ and solve \eqref{eqn:geneval} to obtain $(\tilde{\lambda}_i, e_i), ~i= 1, \ldots, m$
    \STATE Use $\delta$ and $(\tilde{\lambda}_i, e_i)$ to compute the $m-$truncated version of the inner product, denoted $\langle \cdot, \cdot \rangle_{\tilde{K},m}$, via truncating \eqref{eqn:nuggIP}
    \STATE Use $\langle \cdot, \cdot \rangle_{\tilde{K},m}$ to construct the pkGP mean and kernel functions in \eqref{eqn:nuggFnArray}.
    \RETURN $\tilde{\mu}_0(\cdot), ~\tilde{k}_0(\cdot, \cdot)$
    \end{algorithmic}}
\end{algorithm}

\subsection{Boundary Conditions}
{We take as our function of interest the ``non-polynomial function" studied by \cite{LimEtAl02}, so named because it closely resembles a multivariate polynomial. It is defined as
\begin{equation}\label{eqn:nonpoly}
    f(x, y) = \frac{1}{6}\left[ (20 + 5x\sin(5x))(4 + \exp(-5y))-100\right], ~~(x,y) \in T,
\end{equation}
where $T = [0, 1]^2$. Suppose the function is entirely known on the boundary, $T_0 = \partial T$, and we wish to interpolate the function elsewhere. As training data, we take $N = 20$ observations of $f$ in the interior of the domain, chosen by random Latin hypercube design \cite{lhsr}. Interpolators we compare are the projected kernel Gaussian process (pkGP) proposed in this article and the ordinary kriging interpolator. In the absence of a formally-defined projected kernel, in practice one might simply take the ordinary GP predictor and augment the training data with a finite number of pseudo-observations along the boundary, which are available {\em a priori} without having to evaluate $f$. That is to say, whereas the classical kriging interpolator would use training data $\boldsymbol{y}$ obtained via (e.g.) Latin hypercube sampling of the interior of the domain, one might implement ``pseudo-kriging," which is the same as ordinary kriging but with additional training data, $\boldsymbol{y}^* = (\boldsymbol{y}^\top, \boldsymbol{y}_p^\top)^\top$, where $\boldsymbol{y}_p$ are the additional pseudo-observations on the known subset. 

Our proposed pkGP uses an order $M = 15$ Rayleigh-Ritz approximation of the eigenvalues of $k_0$ as found via solving \eqref{eqn:geneval} with Legendre polynomials $\{\xi_j\}_{j=1}^M$. As in Section \ref{subsec:verify}, we build the pseudo-kriging predictor by augmenting the $N$ interior points with $4M = 60$ evenly spaced function values along $T_0$. We use the Mat\'{e}rn kernel with smoothness $\nu = 3/2$ and length-scale $\rho = 1$ for pkGP (prior to projection) and the kriging interpolators; i.e., $k(\boldsymbol{x}, \boldsymbol{y}) = (1 + \|\boldsymbol{x} - \boldsymbol{y}\|\sqrt{3})\exp(-\|\boldsymbol{x} - \boldsymbol{y}\|\sqrt{3})$. The test points at which we evaluate the predictive fidelities of the three are taken to be 81 evenly spaced points in $\{(0.9x, 0.9y) : (x,y) \in \partial T\}$ and 81 evenly spaced points in $\{(0.5x, 0.5y) : (x,y) \in \partial T\}$, so that we are testing near the boundary, and further toward the interior of the domain.}

{Figure \ref{fig:boundaryPerspPlot} plots the test function \eqref{eqn:nonpoly}, along with predicted output from each of the GP interpolators. For further exposition, we plot in Figure \ref{fig:boundaryPreds} the true function evaluations against the predicted values for each GP model. In terms of the root mean squared error (RMSE), ordinary kriging using no boundary information is clearly the worse performer ($RMSE = 0.3544$), as expected. The other two are competitive with each other compared to ordinary kriging, though our proposed pkGP performs the best ($RMSE_{pseudo} = 0.1017$, $RMSE_{pkGP} = 0.0995$).}
\begin{figure}[tb]
    \centering
    \includegraphics[scale= 0.2]{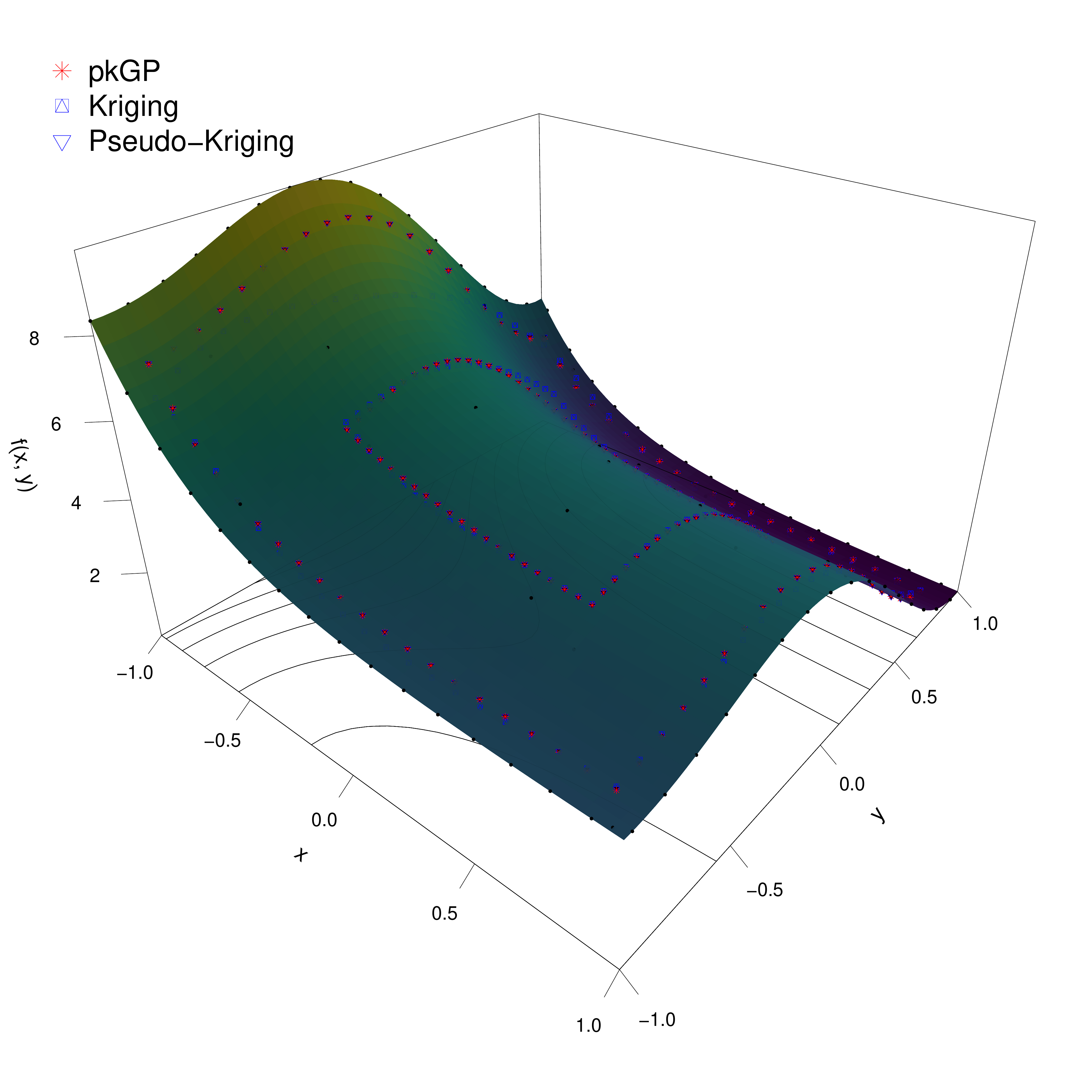}
    \caption{{\em Plot of the test function along with the predicted points from each GP interpolator considered in the known boundary example. The black dots indicate both the training points and the pseudo-observations used for pseudo-kriging.}}
    \label{fig:boundaryPerspPlot}
\end{figure}
\begin{figure}[tb]
    \centering
    \includegraphics[scale= 0.125]{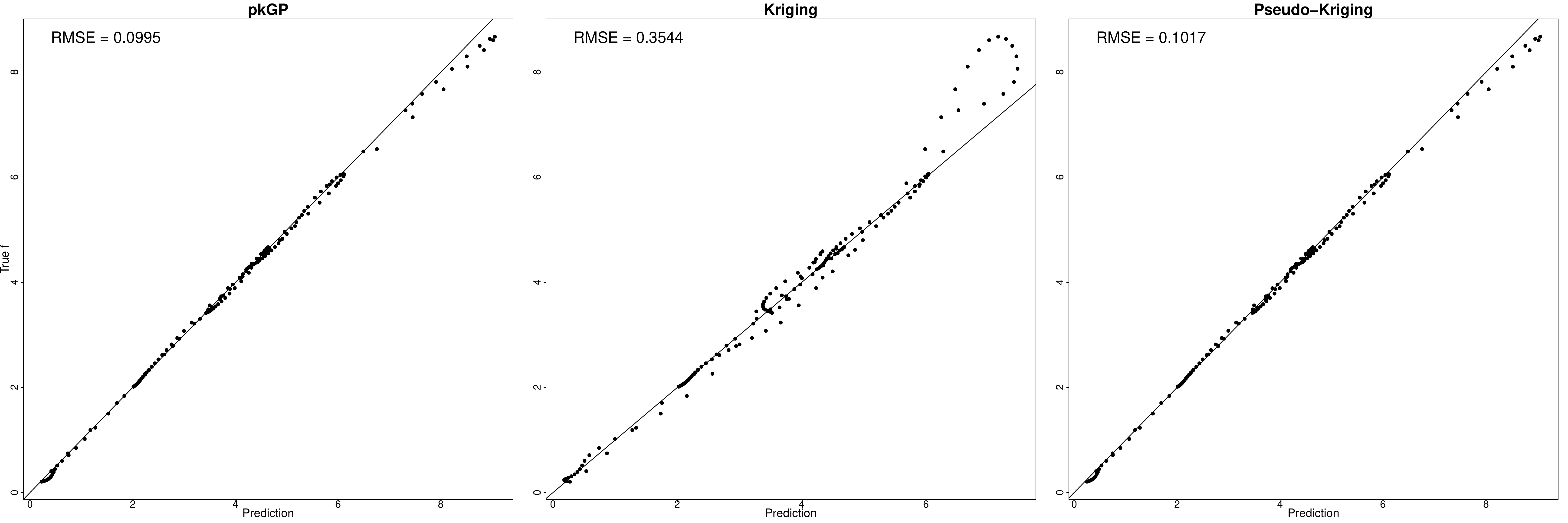}
    \caption{{\em True function values versus predictions for each of the GP interpolators in the known boundary example.}}
    \label{fig:boundaryPreds}
\end{figure}

{To further compare our proposed pkGP to both kriging versions over the entirety of function surfaces, we consider an additional three test functions commonly used in the literature \cite{surjanovic2016virtual, ding2019bdrygp}. These functions are called the ``corner peak" function, the ``product peak" function, and the Rosenbrock function. The functions, denoted $f_{corn}, f_{prod}, f_{rosen},$ respectively, are given by the following:
\begin{eqnarray*}
    f_{corn}(x_1, x_2) &=& \left(1 + \frac{x_1 + x_2}{2}\right)^{-3}, ~x_1, x_2 \in [0, 1]\\
    f_{prod}(x_1, x_2) &=& \prod_{i=1}^2(1 + 10(x_i - 0.25)^2)^{-1}, ~x_1, x_2 \in [0, 1]\\
    f_{rosen}(x_1, x_2) &=& 100(x_2 - x_1^2)^2 + (1-x_1)^2, ~x_1, x_2 \in [0, 1]\\
\end{eqnarray*}
They are plotted in Figure \ref{fig:threeFns}. We use the same Mat\'{e}rn kernel and same number of pseudo-observations for the pseudo-kriging predictor along the boundary. The number of (interior) training points are varied from 10 to 200, where each sample is obtained via Latin hypercube sampling on $[0, 1]^2$.}
\begin{figure}[tb]
    \centering
    \includegraphics[scale= 0.1]{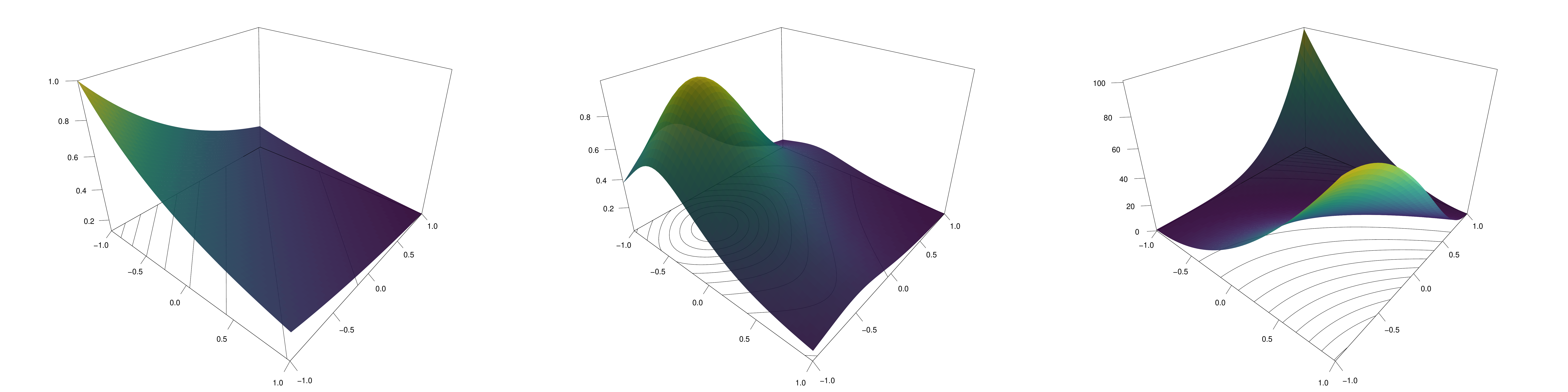}
    \caption{{\em Plot of the three test functions used in the boundary condition example. The left, center, and right panels display the corner peak, product peak, and Rosenbrock functions, respectively.}}
    \label{fig:threeFns}
\end{figure}

{Figure \ref{fig:threeRelErrs} displays the approximate relative errors with respect to the supremum norm. For the corner peak and Rosenbock functions, the proposed pkGP outperforms both ordinary and pseudo-kriging, with the differences generally becoming more pronounced as the training data size increases. Pseudo-kriging and pkGP have indistinguishable relative errors for the product peak function. This behavior is presumably because the product peak function is more well behaved along the boundary compared to the other two. In other words, it appears that the difference between projecting directly onto the subspace versus using a finite number of pseudo observations is more pronounced when the sharply approaches a maximum near its boundary.}\\
\begin{figure}[tb]
    \centering
    \includegraphics[scale= 0.1]{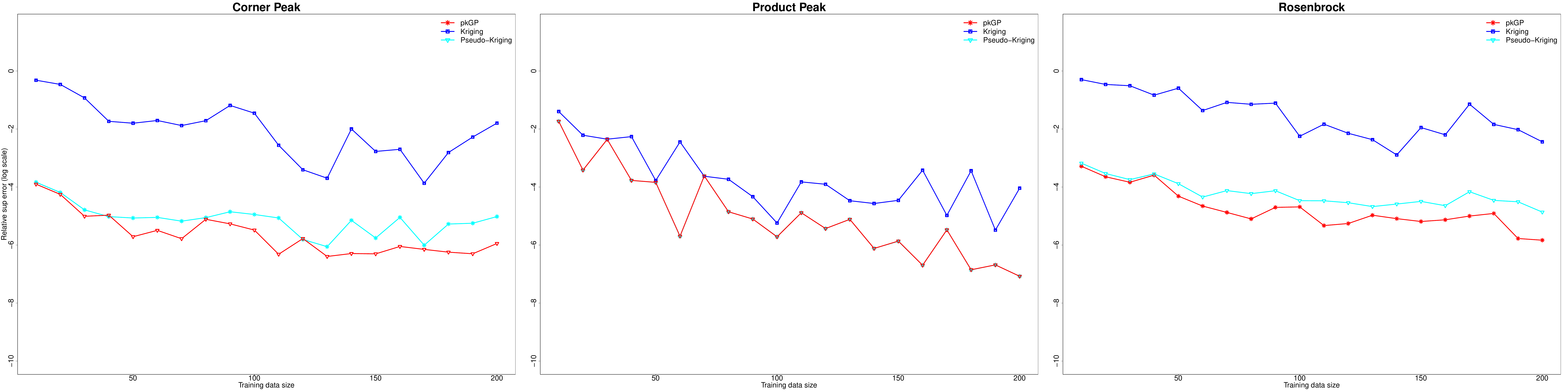}
    \caption{{\em Relative $\|\cdot\|_{\infty}$ errors of each interpolator as a function of training data size. The test functions are displayed in Figure \ref{fig:threeFns}. The relative errors of pkGP and pseudo-kriging are indistinguishable for the product peak function (middle panel).}}
    \label{fig:threeRelErrs}
\end{figure}
{The results from the known boundary illustrations suggest that pkGP and kriging augmented with pseudo-observations may perform similar to each other in certain settings. It is worth emphasizing, though, that pkGP does not use pseudo-training data that would otherwise increase the size of the matrix to be factored and inverted for ordinary kriging. We further compare pkGP to pseudo-kriging in the next example.}

\subsection{Diagonal Conditions}

The projected kernel Gaussian process (pkGP) proposed in this work is not limited to cases of boundary constraints; i.e., $T_0$ is not limited to the boundary, but can be any compact subset of the domain $T$. In this example, we again assume that $T=[-1, 1]^2$. The target function of interest is given by
$$
f(x, y) = y\sqrt{1+x}\cos(\pi y)\sin\Big(\frac{\pi(x-y)}{2}+1\Big) e^{.5 (x + y) ^ 2}, ~~(x, y) \in [-1, 1]^2.
$$

Rather than known boundary conditions, we assume that $f$ is known along the diagonal of $T$, $T_0 = \{(t, t): t \in [-1,1]\} \subset T$. Our goal is to approximate as well as possible the function values along the test set of points near the boundary; i.e., test values contained in  the set $\{(t, t\pm .1):t \in [-0.9, 0.9]\}$. The $N= 20$ training points are obtained via Latin hypercube sampling. We consider pkGP and pseudo-kriging as in the boundary example. The projected kernel is approximated the same as in the previous example with $M= 16$ basis Legendre polynomials for the Rayleigh-Ritz method. We implement also the analogous pseudo-kriging predictor augmented with 16 pseudo-observations along the diagonal.

{Figure \ref{fig:diagplot} plots the target function along with predicted output from both the proposed pkGP and pseudo-kriging. Also plotted are the training points and pseudo-observations. We emphasize that while the pseudo-observations are used with the typical kriging, they are {\em not} used for pkGP, which is already projected onto $T_0$. This suggests computational savings that may be realized via our suggested approach of projecting in function space prior to training. The quality of the approximation is also displayed and quantified in Figure \ref{fig:diagPreds}. In terms of RMSE, pkGP produces substantially more faithful predictions than the pseudo-kriging predictor at the test points ($RMSE_{pkGP} = 0.1413$, $RMSE_{pseudo} = 0.2240$).}
\begin{figure}[tb]
    \centering
    \includegraphics[scale= 0.2]{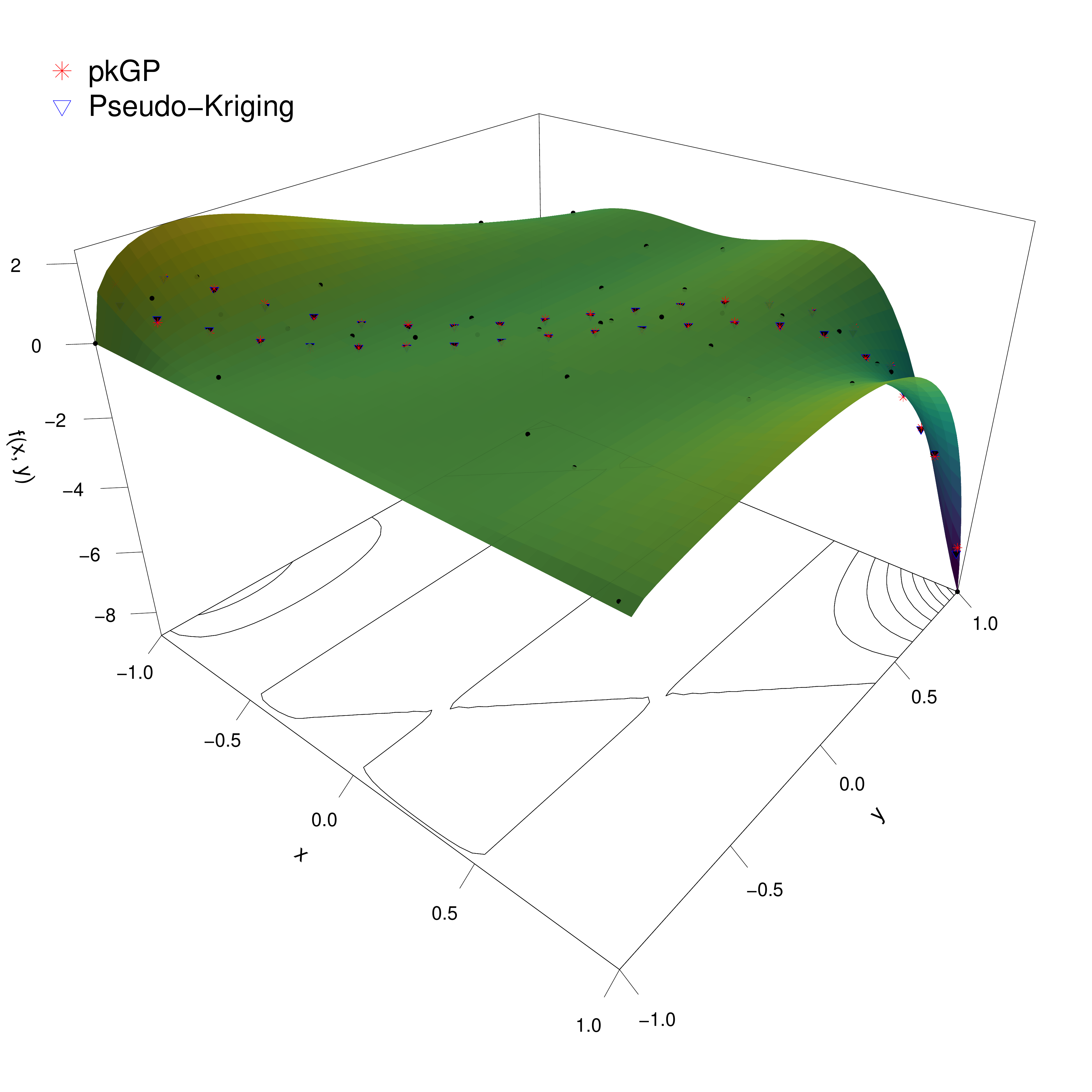}
    \caption{{\em Plot of the test function along with the predicted points from each GP interpolator considered in the known diagonal example. The black dots indicate both the training points and the pseudo-observations used for pseudo-kriging.}}
     \label{fig:diagplot}
\end{figure}
\begin{figure}[tb]
    \centering
    \includegraphics[scale= 0.15]{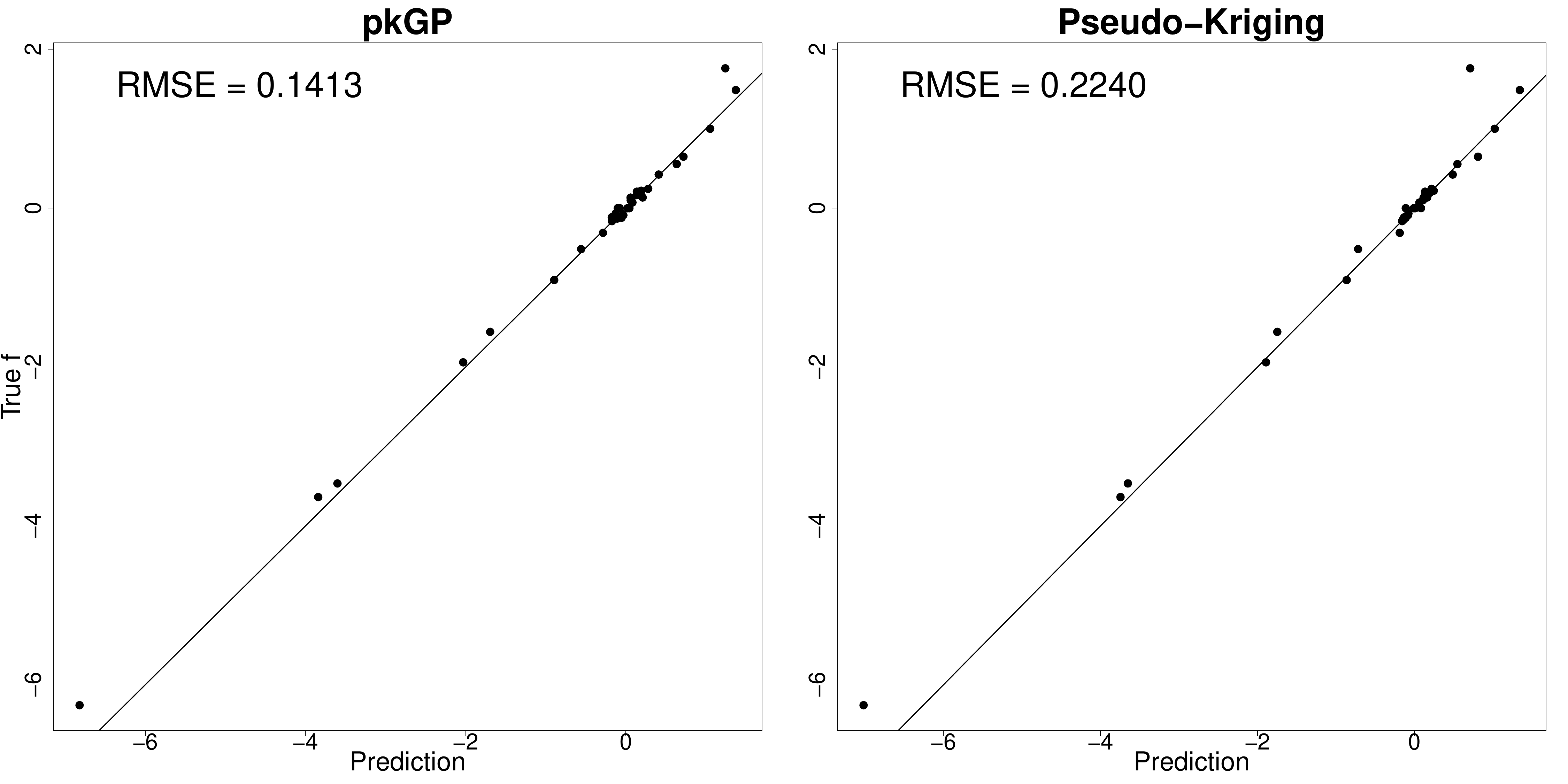}
    \caption{{\em True function values versus predictions for each of the GP interpolators in the known diagonal example.}}
    \label{fig:diagPreds}
\end{figure}

{Since the function is known along the entire uncountable diagonal subset of the domain, we examine the effect of increasing the number of pseudo-observations along the diagonal, as one might do in practice to approximate the infinitely-many known points. Again, such pseudo-observations are not needed (and are in fact redundant) for our proposed pkGP. The known subset is automatically incorporated into the pkGP kernel function via orthogonal projection, so we need not consider increasing the its number of training points. Figure \ref{fig:diagPseudo} plots the relative approximation errors with respect to the $L^2$ norm, $\|\hat{f} - f\|_2/\|f\|_2$, against the number of pseudo-observations for pseudo-kriging. As expected, we see the pseudo-kriging approximation improving as it is conditioned upon more information along the boundary, approaching that of pkGP. However, it does not attain the lower error from the proposed pkGP. This plot in particular is indicative of the result shown in this paper that the projected kernel GP is the limit of finite-dimensional conditioned GPs.}
\begin{figure}[tb]
    \centering
    \includegraphics[scale= 0.1]{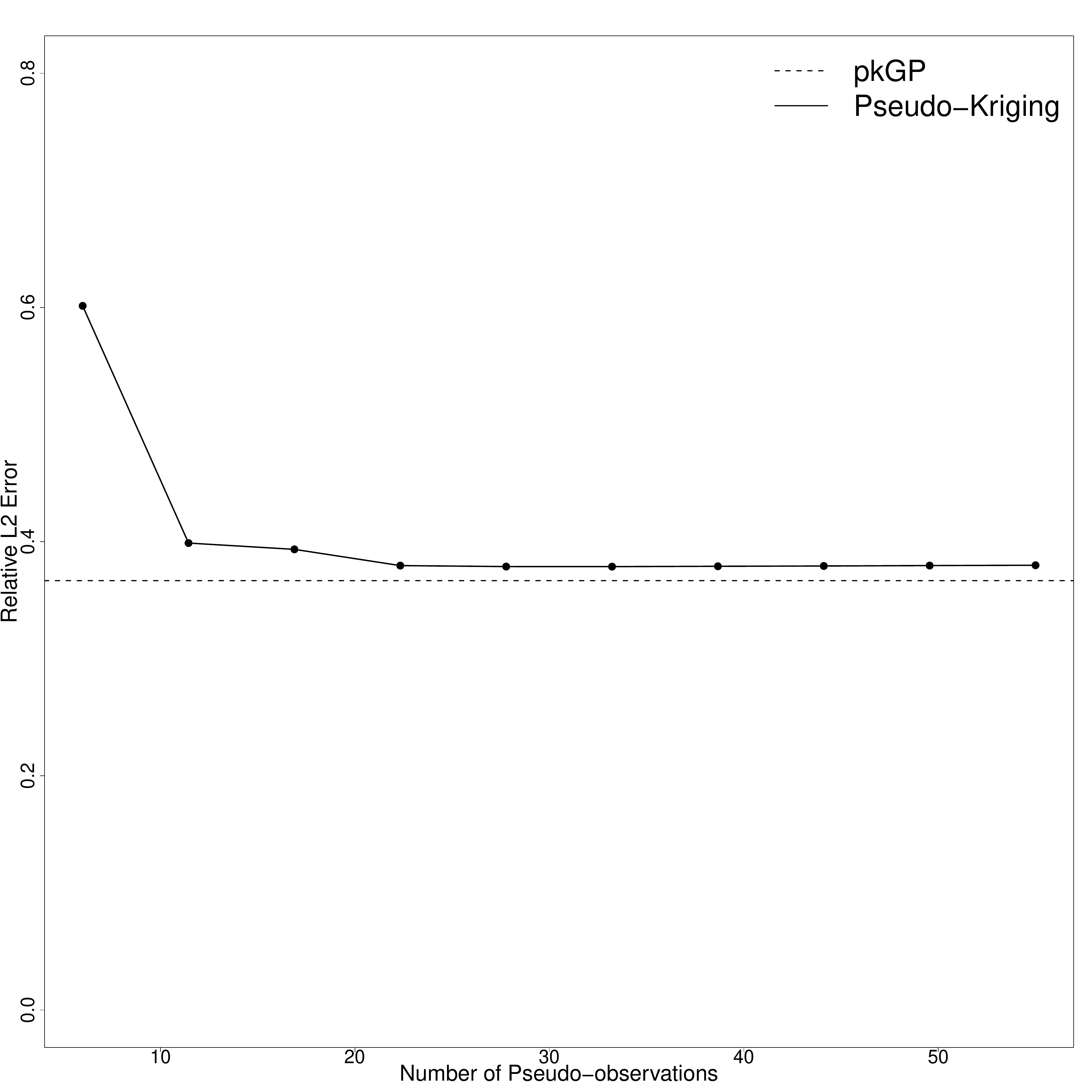}
    \caption{{\em Relative $L^2$ errors of ordinary kriging function approximation in the known diagonal example, as a function of number of pseudo-observations used. The horizontal dashed line is the relative error of pkGP, which does not use any pseudo-observations.}}
    \label{fig:diagPseudo}
\end{figure}

\subsection{Summary}
{Our numerical experiments illustrate the flexibility of our proposed projected kernel Gaussian process to different types of known conditions, boundary constraints and more general subset constraints. In the case of boundary constraints, pkGP outperforms or is otherwise competitive with both ordinary kriging and pseudo-kriging. In the diagonal example where boundary constraints are not appropriate, pkGP still outperforms pseudo-kriging augmented with an increasing number of pseudo-observations along the diagonal. This latter illustration suggests that working with the projected kernel GP directly in function space can yield improved predictions without worrying about the number and locations of pseudo observations nor the associated increase in computational burden. Conversely, the {\em ad hoc} practice of using pseudo-observations to ``boost" ordinary kriging can be theoretically justified as approximating a well-defined (and well behaved) infinite-dimensional process.}

\section{Conclusions and Future Directions}\label{sec:final}

The goal of this paper is to construct and study Gaussian processes which are capable of using information from fairly arbitrary subsets of the domain while requiring minimal assumptions. Using the geometry of orthogonal projections in reproducing kernel Hilbert space, we explicitly define the conditional mean and covariance of Gaussian processes, prove that such processes exist, and that they can be expressed as limits of kriging interpolators with an increasing number of pseudo-observations. Numerical examples illustrate the flexibility of our proposed approach, including its ability to outperform existing alternatives. {Future work in this area might include characterizing the functions/scenarios in which such known information is truly beneficial to incorporate, and accounting for possibly discontinuous functions. There remains also the question of additional regularity conditions that could possibly establish rates of convergence, or even non-asymptotic convergence guarantees.}

\section*{Acknowledgments}
We would like to thank the Editorial Board and two anonymous referees for their time and expertise in reviewing this paper, including bringing to our attention additional references of which we were unaware. The comments and suggestions greatly improved the manuscript.\\

\appendix 
\section{Additional Proofs} 
\subsection{Theorem \ref{thm:iso_thm}} \label{app:thm_3.3}

\begin{proof}
Define $\psi:\mbox{Span}(\{k_s;\, s \in T_0\}) \to \mathcal{H}(T_0)$ by $f \mapsto f|_{T_0}$, which is well-defined and linear. 
Note that for arbitrary $n \ge 1$, $\{t_1,\hdots,t_n\} \subset T_0$, $f = \sum_{i=1}^n a_i k_{t_i}$, and we have
$$\langle k_{t_i},k_{t_j} \rangle_{\mathcal{H}(T)} = k_{t_j}(t_i) = (\psi k_{t_j})(t_i) = \langle \psi k_{t_i}, \psi k_{t_j} \rangle_{\mathcal{H}(T_0)}.$$ 

Using this property, it follows that for $n \ge 1$, $\{t_1,\hdots,t_n\} \subset T_0$, $f = \sum_{i=1}^n a_i k_{t_i}$, and we have
\begin{align*}
\langle f,f \rangle_{\mathcal{H}(T)} & = \sum_{i=1}^n \sum_{j=1}^n a_i a_j \langle k_{t_i}, k_{t_j} \rangle_{\mathcal{H}(T)}\\
&= \sum_{i=1}^n \sum_{j=1}^n a_i a_j \langle \psi k_{t_i}, \psi k_{t_j} \rangle_{\mathcal{H}(T_0)} \\
& = \bigg \langle \psi \Big(\sum_{i=1}^n a_i k_{t_i}\Big), \psi\Big(\sum_{j=1}^n a_j k_{t_j}\Big) \bigg \rangle_{\mathcal{H}(T_0)} = \langle \psi f, \psi f \rangle_{\mathcal{H}(T_0)}.
\end{align*}
Therefore $\psi$ is an isometry. Now, define $\Tilde{\psi}: \mathcal{H}_0^\perp \to \mathcal{H}(T_0)$ by $f \mapsto f|_{T_0}$. By Proposition \ref{prop: perp}, each element of $\mathcal{H}_0^\perp$ when restricted to $T_0$ identifies with an element of $\mathcal{H}(T_0)$. Hence, $\Tilde{\psi}$ is well defined. Again, from Proposition \ref{prop: perp}, it remains to show that the isometry property of $\psi$ can be extended to the closure of $\mbox{Span}(\{k_s;\, s \in T_0\})$, and that there is a one-to-one correspondence between $\mathcal{H}_0^\perp$ and $\mathcal{H}(T_0)$ via $\Tilde{\psi}$. 

Note that one may define any element $f \in \mathcal{H}_0^\perp$ as the limit of some Cauchy sequence $\{f_n\} \subset \mbox{Span}(\{k_s;\, s \in T_0\})$. By the continuity of norms, $\langle f, f \rangle_{\mathcal{H}(T)} = \lim_{n \to \infty} \langle f_n, f_n \rangle_{\mathcal{H}(T)}$. Then, by the isometry property of $\psi$, $\langle f_n, f_n \rangle_{\mathcal{H}(T)} = \langle \psi f_n , \psi f_n  \rangle_{\mathcal{H}(T_0)}$. Since $\Tilde{\psi} = \psi$ on $\mbox{Span}(\{k_s;\, s \in T_0\})$, we have $\langle \psi f_n , \psi f_n  \rangle_{\mathcal{H}(T_0)} = \langle \Tilde{\psi} f_n, \Tilde{\psi} f_n\rangle_{\mathcal{H}(T_0)}$. Therefore,
\begin{align*}
\Big| \langle f, f \rangle_{\mathcal{H}(T)} - \langle \Tilde{\psi}f, \Tilde{\psi}f \rangle_{\mathcal{H}(T_0)} \Big| & = \lim_{n \to \infty} \Big| \langle \Tilde{\psi} f_n , \Tilde{\psi} f_n  \rangle_{\mathcal{H}(T_0)} - \langle \Tilde{\psi}f, \Tilde{\psi}f \rangle_{\mathcal{H}(T_0)}\Big| \\
 & \le \lim_{n \to \infty} \Big| \langle \Tilde{\psi}(f - f_n), \Tilde{\psi}f \rangle_{\mathcal{H}(T_0)}\Big| + \Big| \langle \Tilde{\psi}(f - f_n), \Tilde{\psi}f_n \rangle_{\mathcal{H}(T_0)}\Big| \\
 & \le 2 ||\Tilde{\psi}||^2 \sup_{h \in \{f\} \cup \{f_n\}} ||h|| \lim_{n \to \infty} ||f - f_n||_{\mathcal{H}(T)} = 0.
\end{align*}
Hence, $\Tilde{\psi}$ is an isometry. It remains to show $\Tilde{\psi}$ is one-to-one and onto.

{For $f \in \mathcal{H}_0^\perp$,} $\Tilde{\psi}f \equiv 0$ implies that $f|_{T_0} \equiv 0$ so that $f \in \mathcal{H}_0$. Thus, $f \equiv 0$ {and $\Tilde{\psi}$ is one-to-one.}

{To show that $\tilde{\psi}$ is surjective}, suppose $h \in \mathcal{H}(T_0)$. Then there exists a Cauchy sequence $\{h_n\} \subset \mbox{Span}(\{k_s|_{T_0}; \, s \in T_0\})$ which converges to $h$. {Then one may find a sequence, say $\{f_n\} \subset \mbox{Span}(\{k_s ; s \in T_0\}) \subset \mathcal{H}_0^\perp$, so that $\Tilde{\psi}f_n=h_n$ for all $n$}. Since $\Tilde{\psi}$ is an isometry {and $\{h_n\}$ is Cauchy}, $\{f_n\}$ is Cauchy and therefore has a limit $f \in \mathcal{H}_0^\perp$. It follows that
$$
\Tilde{\psi}f = \Tilde{\psi}\Big(\lim_n f_n\Big) = \lim_n \, \Tilde{\psi} f_n = \lim_n h_n = h.
$$
Thus, $\Tilde{\psi}$ is onto and hence an isomorphism.
\end{proof}

\subsection{Lemma \ref{lem:continuous}} \label{app:lem_3.4}
\begin{proof}
We will use the Kolmogorov-Chentsov theorem \cite[Theorem 2.23]{kall} which states that $\mathbb{X}$ has a continuous version on $\mathbb{R}^d$ taking on values in a complete metric space $(S,\rho)$ if there exists $a,b>0$ such that
$$
E[\rho(X_s,X_t)^a] \le {c}|s-t|^{d+b}, \, s,t \in \mathbb{R}^d,
$$
{for some constant $c$}. Assume {without loss of generality} that $\mathbb{X}$ has zero mean. Define $\rho$ to be the Euclidean norm on $\mathbb{R}$, and recall that for any zero mean Gaussian random variable $Z$ and any even integer $a$,
$$
E[Z^{a}]= C_a E[Z^2]^{a/2},
$$
where $C_a=\prod_{i=1}^{a/2} (2i-1)$. Defining $a$ to be the smallest even integer strictly larger than $\frac{2d}{\gamma}$, we see for any $s,t \in \mathbb{R}^d$,
\begin{align*}
    E[\rho(X_t,X_s)^a]&=E[(X_t-X_s)^a]=C_{a} E[(X_t-X_s)^2]^{a/2}=C_a \big [k(t,t)-2k(t,s)+k(s,s) \big]^{a/2} \\
    & \le C_a|s-t|^{\gamma a/2}=C_a|s-t|^{d+(\gamma a/2-d)}.
\end{align*}
Thus, selecting $b=\gamma a/2-d$, and scaling $\rho$ appropriately, we get the result for a zero mean process. Lastly, the non-zero mean process can be achieved by translating the process by the mean, repeating the procedure above, and noting that the sum of continuous functions is continuous.
\end{proof}

\subsection{Lemma \ref{lem:tight}} \label{app:thm_3.5}
\begin{proof}
Recall the remark in Section 3 in which the mean and covariance of $\mathbb{X}^N$, denoted $\mu^N$ and $k^N$, can be defined as
\begin{align*}
    \mu^N(s) &=\mu(s)+\langle Q_N{\hat{k}_s},Q_N{(\hat{g}-\hat{\mu}}) \rangle_{\mathcal{H}(T_0)}, \\
    k^N_0(s,t) &=k(s,t)-\langle Q_N {\hat{k}_s}, Q_N {\hat{k}_t}\rangle_{\mathcal{H}(T_0)}.
\end{align*}
{We next observe the relationships between the spaces $\mathcal{H}(T), ~\mathcal{H}(T_0),$ and $\mbox{Span}(\{\hat{k}_{t_1}, \ldots, \hat{k}_{t_N}\})$, which can be seen in the following diagram:}

\vspace*{4 mm}
\begin{center}
\begin{tikzcd}
\mathcal{H}(T) \arrow[r, "Q"] & \mathcal{H}_0^{\perp} \arrow[r,"\tilde{\psi}"]
 & \mathcal{H}(T_0) \arrow[r, "Q_N"] & \mbox{Span$(\{ \hat{k}_{t_1},\ldots,\hat{k}_{t_N}\})$}
\end{tikzcd}
\end{center}
\vspace*{4mm}
{This illustrates the fact that, for $f \in \mathcal{H}(T)$, $f_N = (Q_N \circ \tilde{\psi} \circ Q)f$. From this we can see that $\|Q_N\hat{f}\|_{\mathcal{H}(T)} = \|f_N\|_{\mathcal{H}(T_0)} \leq \|f\|_{\mathcal{H}(T)}$.} Now, {making use of Lemma \ref{lemma:kHolder}}, observe that for $s_0 \in T$,
\begin{align*}
|k^N_0(s_0,s) - k^N_0(s_0,t)| &\le |k(s_0,s)-k(s_0,t)| + |\langle Q_N {\hat{k}_{s_0}}, Q_N {(\hat{k}_s-\hat{k}_t)}\rangle_{\mathcal{H}(T_0)}| \\
&\le C|s-t|^{\gamma} + ||Q_N {\hat{k}_{s_0}} ||_{\mathcal{H}(T_0)} || Q_N({\hat{k}_s-\hat{k}_t)}||_{\mathcal{H}(T_0)} \\
&\le C|s-t|^{\gamma} + ||{\hat{k}_{s_0}}||_{\mathcal{H}(T_0)} ||k_s-k_t||_{\mathcal{H}(T)} \\
& \le C|s-t|^{\gamma} + C'|s-t|^{\gamma/2}\\
&\le \Tilde{C} |s-t|^{\gamma/2},
\end{align*}
where the first inequality follows from the triangle inequality, the final inequality follows form the boundedness of $T$, and $\Tilde{C}$ does not depend on $s_0$ or $N$.  Since $k$ itself is $\gamma-$H\"older continuous, it follows that $k^N_0$ is $\gamma/2-$H\"older continuous on $T \times T$ uniformly in $N$. Furthermore, $\mu^N \to \mu$ uniformly where we again use the fact that {$K$} is uniformly $\gamma/2-$H\"older continuous on $\{Q_N(g-\mu)\}_{N=1}^\infty$. Therefore, $\{\mathbb{X}^N\}_{N=1}^\infty$ is tight.
\end{proof}

\bibliographystyle{plain}
\bibliography{references-final-ar}

\end{document}